\documentclass[12pt,english]{amsart}
\usepackage[latin1]{inputenc}
\usepackage{a4wide}
\usepackage{graphicx}
\usepackage{amssymb}
\usepackage{babel}
\usepackage{psfrag}

 \theoremstyle{plain}    
 \newtheorem{thm}{Theorem}[section]
 \numberwithin{equation}{section}
 \numberwithin{figure}{section}
 \theoremstyle{plain}
 \theoremstyle{plain}    
 \newtheorem{conjecture}[thm]{Conjecture}
 \theoremstyle{plain}    
 \newtheorem{lem}[thm]{Lemma}
 \theoremstyle{plain}    
 \newtheorem{cor}[thm]{Corollary}
 \theoremstyle{remark}    
 \newtheorem*{remark}{Remark}

\def\la{\lambda}
\def\si{\sigma}
\def\al{\alpha}
\def\be{\beta}
\def\et{\eta}
\def\P{\mathcal P}
\def\sgn{\operatorname{sgn}}
\newcommand{\DC}{\Delta^{(\varLambda-\mathcal{V})}}
\newcommand{\DO}{\Delta^{(\mathcal{V}-\varLambda)}}

\makeatletter
\newfont{\footsc}{cmcsc10 at 8truept}
\newfont{\footbf}{cmbx10 at 8truept}
\newfont{\footrm}{cmr10 at 10truept}

\renewcommand{\ps@plain}{%
\renewcommand{\@oddfoot}{\footsc the electronic journal of combinatorics
  {\footbf 11}(2) (2005), \#R16\hfil\footrm\thepage}
\renewcommand{\@evenfoot}{\footsc the electronic journal of combinatorics
  {\footbf 11}(2) (2005), \#R16\hfil\footrm\thepage}}

\makeatother

\catcode`\@=11

\thinlines
\newskip\Einheit \Einheit=0.6cm
\newcount\xcoord \newcount\ycoord
\newdimen\xdim \newdimen\ydim \newdimen\PfadD@cke \newdimen\Pfadd@cke
\def\PfadDicke#1{\PfadD@cke#1 \divide\PfadD@cke by2 \Pfadd@cke\PfadD@cke \multiply\PfadD@cke by2}
\long\def\LOOP#1\REPEAT{\def\BODY{#1}\ITERATE}
\def\ITERATE{\BODY \let\next\ITERATE \else\let\next\relax\fi \next}
\let\REPEAT=\fi
\def\Punkt{\hbox{\raise-2pt\hbox to0pt{\hss\scriptsize$\bullet$\hss}}}
\def\DuennPunkt(#1,#2){\unskip
  \raise#2 \Einheit\hbox to0pt{\hskip#1 \Einheit
          \raise-2.5pt\hbox to0pt{\hss\normalsize$\bullet$\hss}\hss}}
\def\NormalPunkt(#1,#2){\unskip
  \raise#2 \Einheit\hbox to0pt{\hskip#1 \Einheit
          \raise-3pt\hbox to0pt{\hss\large$\bullet$\hss}\hss}}
\def\DickPunkt(#1,#2){\unskip
  \raise#2 \Einheit\hbox to0pt{\hskip#1 \Einheit
          \raise-4pt\hbox to0pt{\hss\Large$\bullet$\hss}\hss}}
\def\Kreis(#1,#2){\unskip
  \raise#2 \Einheit\hbox to0pt{\hskip#1 \Einheit
          \raise-4pt\hbox to0pt{\hss\Large$\circ$\hss}\hss}}
\def\Diagonale(#1,#2)#3{\unskip\leavevmode
  \xcoord#1\relax \ycoord#2\relax
      \raise\ycoord \Einheit\hbox to0pt{\hskip\xcoord \Einheit
         \unitlength\Einheit
         \line(1,1){#3}\hss}}
\def\AntiDiagonale(#1,#2)#3{\unskip\leavevmode
  \xcoord#1\relax \ycoord#2\relax 
      \raise\ycoord \Einheit\hbox to0pt{\hskip\xcoord \Einheit
         \unitlength\Einheit
         \line(1,-1){#3}\hss}}
\def\Pfad(#1,#2),#3\endPfad{\unskip\leavevmode
  \xcoord#1 \ycoord#2 \thicklines\ZeichnePfad#3\endPfad\thinlines}
\def\ZeichnePfad#1{\ifx#1\endPfad\let\next\relax
  \else\let\next\ZeichnePfad
    \ifnum#1=1
      \raise\ycoord \Einheit\hbox to0pt{\hskip\xcoord \Einheit
         \vrule height\Pfadd@cke width1 \Einheit depth\Pfadd@cke\hss}%
      \advance\xcoord by 1
    \else\ifnum#1=2
      \raise\ycoord \Einheit\hbox to0pt{\hskip\xcoord \Einheit
        \hbox{\hskip-\PfadD@cke\vrule height1 \Einheit width\PfadD@cke depth0pt}\hss}%
      \advance\ycoord by 1
    \else\ifnum#1=3
      \raise\ycoord \Einheit\hbox to0pt{\hskip\xcoord \Einheit
         \unitlength\Einheit
         \line(1,1){1}\hss}
      \advance\xcoord by 1
      \advance\ycoord by 1
    \else\ifnum#1=4
      \raise\ycoord \Einheit\hbox to0pt{\hskip\xcoord \Einheit
         \unitlength\Einheit
         \line(1,-1){1}\hss}
      \advance\xcoord by 1
      \advance\ycoord by -1
    \else\ifnum#1=5
      \advance\xcoord by -1
      \raise\ycoord \Einheit\hbox to0pt{\hskip\xcoord \Einheit
         \vrule height\Pfadd@cke width1 \Einheit depth\Pfadd@cke\hss}%
    \else\ifnum#1=6
      \advance\ycoord by -1
      \raise\ycoord \Einheit\hbox to0pt{\hskip\xcoord \Einheit
        \hbox{\hskip-\PfadD@cke\vrule height1 \Einheit width\PfadD@cke depth0pt}\hss}%
    \else\ifnum#1=7
      \advance\xcoord by -1
      \advance\ycoord by -1
      \raise\ycoord \Einheit\hbox to0pt{\hskip\xcoord \Einheit
         \unitlength\Einheit
         \line(1,1){1}\hss}
    \else\ifnum#1=8
      \advance\xcoord by -1
      \advance\ycoord by +1
      \raise\ycoord \Einheit\hbox to0pt{\hskip\xcoord \Einheit
         \unitlength\Einheit
         \line(1,-1){1}\hss}
    \fi\fi\fi\fi
    \fi\fi\fi\fi
  \fi\next}
\def\hSSchritt{\leavevmode\raise-.4pt\hbox to0pt{\hss.\hss}\hskip.2\Einheit
  \raise-.4pt\hbox to0pt{\hss.\hss}\hskip.2\Einheit
  \raise-.4pt\hbox to0pt{\hss.\hss}\hskip.2\Einheit
  \raise-.4pt\hbox to0pt{\hss.\hss}\hskip.2\Einheit
  \raise-.4pt\hbox to0pt{\hss.\hss}\hskip.2\Einheit}
\def\vSSchritt{\vbox{\baselineskip.2\Einheit\lineskiplimit0pt
\hbox{.}\hbox{.}\hbox{.}\hbox{.}\hbox{.}}}
\def\DSSchritt{\leavevmode\raise-.4pt\hbox to0pt{%
  \hbox to0pt{\hss.\hss}\hskip.2\Einheit
  \raise.2\Einheit\hbox to0pt{\hss.\hss}\hskip.2\Einheit
  \raise.4\Einheit\hbox to0pt{\hss.\hss}\hskip.2\Einheit
  \raise.6\Einheit\hbox to0pt{\hss.\hss}\hskip.2\Einheit
  \raise.8\Einheit\hbox to0pt{\hss.\hss}\hss}}
\def\dSSchritt{\leavevmode\raise-.4pt\hbox to0pt{%
  \hbox to0pt{\hss.\hss}\hskip.2\Einheit
  \raise-.2\Einheit\hbox to0pt{\hss.\hss}\hskip.2\Einheit
  \raise-.4\Einheit\hbox to0pt{\hss.\hss}\hskip.2\Einheit
  \raise-.6\Einheit\hbox to0pt{\hss.\hss}\hskip.2\Einheit
  \raise-.8\Einheit\hbox to0pt{\hss.\hss}\hss}}
\def\SPfad(#1,#2),#3\endSPfad{\unskip\leavevmode
  \xcoord#1 \ycoord#2 \ZeichneSPfad#3\endSPfad}
\def\ZeichneSPfad#1{\ifx#1\endSPfad\let\next\relax
  \else\let\next\ZeichneSPfad
    \ifnum#1=1
      \raise\ycoord \Einheit\hbox to0pt{\hskip\xcoord \Einheit
         \hSSchritt\hss}%
      \advance\xcoord by 1
    \else\ifnum#1=2
      \raise\ycoord \Einheit\hbox to0pt{\hskip\xcoord \Einheit
        \hbox{\hskip-2pt \vSSchritt}\hss}%
      \advance\ycoord by 1
    \else\ifnum#1=3
      \raise\ycoord \Einheit\hbox to0pt{\hskip\xcoord \Einheit
         \DSSchritt\hss}
      \advance\xcoord by 1
      \advance\ycoord by 1
    \else\ifnum#1=4
      \raise\ycoord \Einheit\hbox to0pt{\hskip\xcoord \Einheit
         \dSSchritt\hss}
      \advance\xcoord by 1
      \advance\ycoord by -1
    \else\ifnum#1=5
      \advance\xcoord by -1
      \raise\ycoord \Einheit\hbox to0pt{\hskip\xcoord \Einheit
         \hSSchritt\hss}%
    \else\ifnum#1=6
      \advance\ycoord by -1
      \raise\ycoord \Einheit\hbox to0pt{\hskip\xcoord \Einheit
        \hbox{\hskip-2pt \vSSchritt}\hss}%
    \else\ifnum#1=7
      \advance\xcoord by -1
      \advance\ycoord by -1
      \raise\ycoord \Einheit\hbox to0pt{\hskip\xcoord \Einheit
         \DSSchritt\hss}
    \else\ifnum#1=8
      \advance\xcoord by -1
      \advance\ycoord by 1
      \raise\ycoord \Einheit\hbox to0pt{\hskip\xcoord \Einheit
         \dSSchritt\hss}
    \fi\fi\fi\fi
    \fi\fi\fi\fi
  \fi\next}
\def\Koordinatenachsen(#1,#2){\unskip
 \hbox to0pt{\hskip-.5pt\vrule height#2 \Einheit width.5pt depth1 \Einheit}%
 \hbox to0pt{\hskip-1 \Einheit \xcoord#1 \advance\xcoord by1
    \vrule height0.25pt width\xcoord \Einheit depth0.25pt\hss}}
\def\Koordinatenachsen(#1,#2)(#3,#4){\unskip
 \hbox to0pt{\hskip-.5pt \ycoord-#4 \advance\ycoord by1
    \vrule height#2 \Einheit width.5pt depth\ycoord \Einheit}%
 \hbox to0pt{\hskip-1 \Einheit \hskip#3\Einheit 
    \xcoord#1 \advance\xcoord by1 \advance\xcoord by-#3 
    \vrule height0.25pt width\xcoord \Einheit depth0.25pt\hss}}
\def\Gitter(#1,#2){\unskip \xcoord0 \ycoord0 \leavevmode
  \LOOP\ifnum\ycoord<#2
    \loop\ifnum\xcoord<#1
      \raise\ycoord \Einheit\hbox to0pt{\hskip\xcoord \Einheit\Punkt\hss}%
      \advance\xcoord by1
    \repeat
    \xcoord0
    \advance\ycoord by1
  \REPEAT}
\def\Gitter(#1,#2)(#3,#4){\unskip \xcoord#3 \ycoord#4 \leavevmode
  \LOOP\ifnum\ycoord<#2
    \loop\ifnum\xcoord<#1
      \raise\ycoord \Einheit\hbox to0pt{\hskip\xcoord \Einheit\Punkt\hss}%
      \advance\xcoord by1
    \repeat
    \xcoord#3
    \advance\ycoord by1
  \REPEAT}
\def\Label#1#2(#3,#4){\unskip \xdim#3 \Einheit \ydim#4 \Einheit
  \def\lo{\advance\xdim by-.5 \Einheit \advance\ydim by.5 \Einheit}%
  \def\llo{\advance\xdim by-.25cm \advance\ydim by.5 \Einheit}%
  \def\loo{\advance\xdim by-.5 \Einheit \advance\ydim by.25cm}%
  \def\o{\advance\ydim by.25cm}%
  \def\ro{\advance\xdim by.5 \Einheit \advance\ydim by.5 \Einheit}%
  \def\rro{\advance\xdim by.25cm \advance\ydim by.5 \Einheit}%
  \def\roo{\advance\xdim by.5 \Einheit \advance\ydim by.25cm}%
  \def\l{\advance\xdim by-.30cm}%
  \def\r{\advance\xdim by.30cm}%
  \def\lu{\advance\xdim by-.5 \Einheit \advance\ydim by-.6 \Einheit}%
  \def\llu{\advance\xdim by-.25cm \advance\ydim by-.6 \Einheit}%
  \def\luu{\advance\xdim by-.5 \Einheit \advance\ydim by-.30cm}%
  \def\u{\advance\ydim by-.30cm}%
  \def\ru{\advance\xdim by.5 \Einheit \advance\ydim by-.6 \Einheit}%
  \def\rru{\advance\xdim by.25cm \advance\ydim by-.6 \Einheit}%
  \def\ruu{\advance\xdim by.5 \Einheit \advance\ydim by-.30cm}%
  #1\raise\ydim\hbox to0pt{\hskip\xdim
     \vbox to0pt{\vss\hbox to0pt{\hss$#2$\hss}\vss}\hss}%
}
\catcode`\@=12

\begin{document}

\def\figurename{\small Figure}
\newcommand{\set}{\operatorname{set}}

\newbox\Adr
\setbox\Adr\vbox{
\centerline{\sc F. Caselli$^{*\#}$, C. Krattenthaler$^*$, B. Lass$^*$
and P. Nadeau$^\dagger$}
\vskip18pt
\centerline{*Institut Camille Jordan, Universit\'e Claude Bernard
Lyon-I,}
\centerline{21, avenue Claude Bernard, F-69622 Villeurbanne Cedex,
France.}
\centerline{E-mail: \texttt{
(caselli,kratt,lass)@euler.univ-lyon1.fr}}
\vskip18pt
\centerline{$^\dagger$Laboratoire de Recherche en Informatique,
Universit\'e Paris-Sud}
\centerline{91405 Orsay Cedex, France}
\centerline{E-mail: \texttt{ nadeau@lri.fr}}
\vskip18pt
\centerline{\small 
Submitted: Feb 17, 2005; Accepted: Mar 14, 2005; Published: Apr 6, 2005}
}

\author{\box\Adr}


\title{On the number of fully packed loop configurations
with a fixed associated matching}

\address{$^*$Institut Camille Jordan,
Universit\'e Claude Bernard Lyon-I,
21, avenue Claude Ber\-nard,
F-69622 Villeurbanne Cedex, France}
\address{$^\dagger$Laboratoire de Recherche en Informatique,
Universit\'e Paris-Sud, 91405 Orsay Cedex,
France}

\dedicatory{Dedicated to Richard Stanley}

\thanks{$^*$Research supported  by EC's IHRP Programme, grant HPRN-CT-2001-00272, ``Algebraic Combinatorics in Europe''. The second author
was also partially supported by
the ``Algebraic Combinatorics" Programme during Spring 2005 
of the Institut Mittag--Leffler of the Royal Swedish Academy of Sciences.}
\thanks{$^{\#}$Current address: Dipartimento di Matematica,
Universit\`a di Roma La Sapienza, P.le A. Moro 3, I-00185 Roma, Italy}


\begin{abstract}
We show that the number of fully packed loop configurations corresponding
to a matching with $m$ nested arches is polynomial in $m$
if $m$ is large enough, thus essentially
proving two conjectures by Zuber [\emph{Electronic J. Combin.}
\textbf{11}(1) (2004), 
Article~\#R13].
\end{abstract}

\maketitle
\thispagestyle{plain}

\section{Introduction}

In this paper we continue the enumerative study of fully packed loop
configurations corresponding to a prescribed matching begun
by the first two authors in
\cite{caskra}, where we proved two conjectures by Zuber \cite{zuber}
on this subject matter. (See also \cite{degier,difra,difra1,difra2} 
for related results.) 
The interest in this study originates in conjectures by Razumov and
Stroganov \cite{RaStAZ}, and by Mitra, Nienhuis, de Gier and Batchelor
\cite{Mitra}, which predict that the coordinates
of the groundstate vectors of certain Hamiltonians in the dense
$O(1)$ loop model are given by the number of fully packed loop
configurations corresponding to particular matchings. Another motivation
comes from the well-known fact (see e.g.\ \cite[Sec.~3]{degier}) that
fully packed loop configurations are 
in bijection with configurations in the six vertex
model, which, in their turn, are
in bijection with alternating sign matrices, and, thus, the
enumeration of fully packed loop configurations corresponding to a prescribed
matching constitutes an interesting refinement of the
enumeration of configurations in the six vertex model or of
alternating sign matrices.

Here we consider configurations with a growing number
of nested arches. We show that the number of configurations is
polynomial in the number of nested arches, thus proving two further
conjectures of Zuber from \cite{zuber}.

In order to explain these conjectures, we have to briefly recall the relevant
definitions from \cite{caskra,zuber}. 
The \emph{fully packed loop model} (FPL model, for short; see for
example \cite{FPL}) is a model of (not necessarily closed) polygons
on a lattice such that each vertex of the lattice is on exactly one
polygon. Whether or not these polygons are closed, we will 
refer to them as \emph{loops}. 

\begin{figure}[htb]

{\centering \includegraphics{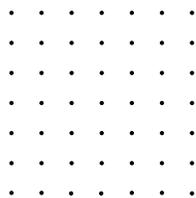} \par}

\caption{\rm The square grid \protect\( Q_{7}\protect \)\label{fpl22}}
\end{figure}  

\begin{figure}[htb] {\centering \includegraphics{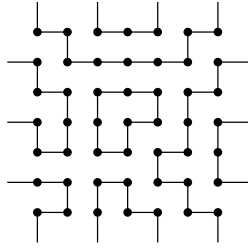} \par}

\caption{An FPL configuration on $Q_7$ with periodic boundary conditions\label{fpl24}} \end{figure}

Throughout this article, we consider
this model on the square grid of side length $n-1$, which we denote
by $Q_{n}$. See Figure~\ref{fpl22} for a picture of $Q_7$. 
The polygons consist of horizontal or vertical edges connecting
vertices of $Q_n$, and edges that lead outside of $Q_n$ from a vertex
along the border of $Q_n$, see Figure~\ref{fpl24} for an example
of an allowed configuration in the FPL model. We call the edges that
stick outside of $Q_n$ {\it external links}. The reader is referred
to Figure~\ref{fpl1} for an illustration of the external links of
the square $Q_{11}$. (The labels should be ignored at this point.)
It should be noted that the four corner points are incident to a horizontal
\emph{and} a vertical external link. We shall be interested here in
allowed configurations in the FPL model, in the sequel referred to
as \emph{FPL configurations}, with \emph{periodic boundary conditions}.
These are FPL configurations where, around the border of $Q_n$, every
other external link of $Q_n$ is part of a polygon. The FPL configuration
in Figure~\ref{fpl24} is in fact a configuration with periodic boundary
conditions. 

\pagenumbering{arabic}
\addtocounter{page}{1}


\medskip

 \medskip

\begin{figure}[htb]

{\centering \includegraphics{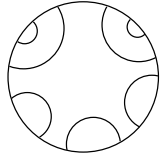} \par}

\caption{The chord diagram representation of a matching\label{fpl25}} \end{figure}

Every FPL configuration defines in a natural way a (non-crossing) 
matching of the
external links by matching those which are on the same polygon (loop).
 We are interested
in the number of FPL configurations corresponding to a fixed matching.
Thanks to a theorem of Wieland \cite{wie} (see Theorem~\ref{thm:wie}),
this number is invariant if the matching is rotated around $Q_n$. This
allows one to represent a matching in form of a chord diagram
of $2n$ points placed around a circle, see
Figure~\ref{fpl25} for the chord diagram representation of the matching
corresponding to the FPL configuration in Figure~\ref{fpl24}.

\begin{figure}[htb]
\psfrag{m}{$m$} \psfrag{X}{$X$}
\begin{center}\includegraphics{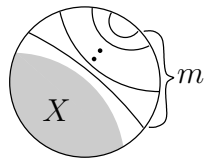}\end{center}
\caption{\label{beau8}The matching composed out of a matching $X$ and
$m$ nested arches}
\end{figure}

The two conjectures by Zuber which we address in this paper concern
FPL configurations corresponding to a matching with $m$ nested
arches. More precisely, let $X$ represent a fixed (non-crossing) 
matching with $n-m$
arches. By adding $m$ nested arches, we obtain a certain matching. 
(See Figure~\ref{beau8} for a schematic
picture of the matching which is composed in this way.) The first of
Zuber's conjectures states that the
number of FPL configurations which has this matching as associated
matching is polynomial
in $m$. In fact, the complete statement is even more precise. 
It makes use of the fact that to any matching $X$ one can associate a
Ferrers diagram $\la(X)$ in a natural way (see Section~\ref{sec:Ferrers} 
for a detailed explanation).

\begin{conjecture}[{\cite[Conj.~6]{zuber}}] \label{conj:zub1}
Let $X$ be a given non-crossing matching with $n-m$ arches, and let $X\cup m$
denote the matching arising from $X$ by adding $m$ nested arches. 
Then the number $A_{X}(m)$ of FPL configurations which have
$X\cup m$ as associated matching is equal to
$\frac{1}{|X|!}P_{X}(m)$,
where $P_{X}(m)$ is a polynomial of degree $\vert \la(X)\vert$ 
with integer coefficients,
and its highest degree coefficient is equal to $\dim(\la(X))$. Here, 
$\vert \la(X)\vert$ denotes the size of the Ferrers diagram $\la(X)$, 
and $\dim(\la(X))$ denotes the dimension of the irreducible
representation of the symmetric group
$S_{\vert \la(X)\vert}$ indexed by the Ferrers diagram $\la(X)$
{\rm(}which is given by the hook formula; see 
\eqref{eq:hook}{\rm)}.
\end{conjecture}

\begin{figure}[htb]
\psfrag{m}{$m$} \psfrag{X}{$X$} \psfrag{Y}{$Y$}
\begin{center}\includegraphics{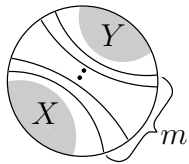}\end{center}
\caption{\label{beau9}Squeezing $m$ nested
arches between two matchings $X$ and $Y$}
\end{figure}

The second conjecture of Zuber generalizes Conjecture~\ref{conj:zub1}
to the case where a bundle of nested arches is squeezed between two
given matchings. More precisely, let $X$ and $Y$ be two given
(non-crossing)
matchings. We produce a new matching by placing $X$ and $Y$ along our
circle that we use for representing matchings, together with $m$
nested arches which we place in between. (See Figure~\ref{beau9} for a
schematic picture.) We denote this matching by $X\cup m\cup Y$.

\begin{conjecture}[{\cite[Conj.~7]{zuber}}] \label{conj:zub2}
Let $X$ and $Y$ be two non-crossing matchings. Then the number
$A_{X,Y}(m)$ of FPL configurations which have
$X\cup m\cup Y$ as associated matching is equal to
$\frac{1}{|\la(X)|!\,|\la(Y)|!}P_{X,Y}(m)$, where $P_{X,Y}(m)$
is a polynomial of degree $|\la(X)|+|\la(Y)|$ with integer coefficients, and
its highest degree coefficient is equal to $\dim(\la(X))\cdot\dim(\la(Y))$.
\end{conjecture}

It is clear that Conjecture~\ref{conj:zub2} 
is a generalization of Conjecture~\ref{conj:zub1},
since $A_{X}(m)=A_{X,\emptyset}(m)$ for any non-crossing matching $X$,
where $\emptyset$ denotes the empty matching. Nevertheless, we
shall treat both conjectures separately, because this will allow us to
obtain, in fact, sharper results
than just the statements in the
conjectures, with our 
result covering Conjecture~\ref{conj:zub1} --- see Theorem~\ref{main}
and Section~\ref{secCon1} --- 
being more precise than the
corresponding result concerning Conjecture~\ref{conj:zub2} --- see
Theorem~\ref{main2a}. We must stress at this point that, while
we succeed to prove Conjecture~\ref{conj:zub1} completely, we are able
to prove Conjecture~\ref{conj:zub2} only for ``large'' $m$, see
the end of Section~\ref{secGen} for the precise statement. There we
also give an explanation of the difficulty of closing the gap.

\medskip
We conclude the introduction by outlining the proofs of our
results, and by explaining the
organisation of our paper at the same time.
All notation
and prerequisites that we are going to use in these proofs
are summarized in Section~\ref{sec:2} below.

Our proofs are based on two observations due to de Gier in
\cite[Sec.~3]{degier} (as are the proofs in
\cite{caskra,difra,difra2}): if one considers the FPL configurations
corresponding to a given matching which has a big number of nested
arches, there are many edges which are occupied by {\it any} such FPL
configuration. We explain this observation, with focus on our
particular problem, in Section~\ref{secfix}. As a consequence, we can
split our enumeration problem into the problem of enumerating
configurations in two separate subregions of $Q_n$, see the
explanations accompanying Figure~\ref{beau15}, respectively
Figure~\ref{beau24}. While one 
of the regions does not depend on $m$, the others grow with $m$.
It remains the task of establishing that the number of configurations
in the latter subregions grows {\it polynomially} with $m$.
In order to do so,
we use the second observation of de Gier, namely the existence of 
a bijection between FPL
configurations (subject to certain constraints on the edges) 
and rhombus tilings, see the proofs of Theorem~\ref{main} and
Lemma~\ref{lem:1}. In the case of
Conjecture~\ref{conj:zub1}, the rhombus tilings can be enumerated
by an application of the hook-content
formula (recalled in Theorem~\ref{thm:rhla}), while in the case of
Conjecture~\ref{conj:zub2} we use a standard correspondence between
rhombus tilings and non-intersecting lattice paths, followed by an
application of the Lindstr\"om--Gessel--Viennot theorem (recalled in
Lemma~\ref{gv}), to obtain
a determinant for the number of rhombus tilings, see the proof of
Lemma~\ref{lem:1}.
In both cases, the polynomial nature of the number of rhombus tilings
is immediately obvious, if $m$ is ``large enough.''
To cover the case of ``small'' $m$ of Conjecture~\ref{conj:zub1} 
as well, we employ a somewhat
indirect argument, which is based on a variation of the above
reasoning, see Section~\ref{secCon1}.
Finally, for the proof of the more specific assertions in
Conjectures~\ref{conj:zub1} and \ref{conj:zub2} on the integrality of
the coefficients of the polynomials (after renormalization) and on the
leading coefficient, we need several technical lemmas
(to be precise, Lemmas~\ref{imposs}, \ref{lem:de} and \ref{lem:3}).
These are implied by Theorem~\ref{th:main}
(see also Corollary~\ref{cor}), which is the subject of
Section~\ref{sec:aux}.


\bigskip

\section{Preliminaries} \label{sec:2}
\subsection{Notation and conventions concerning FPL configurations}

\begin{figure}[htb] 

\psfrag{1}{$1$} 
\psfrag{2}{$2$} 
\psfrag{3}{$3$} 
\psfrag{-1}{$-1$} 
\psfrag{0}{$0$} 
\psfrag{...}{\dots}
\psfrag{-n}{$-n$} 
\psfrag{2n}{$2n$} 
\psfrag{2n-1}{$2n-1$} 
\psfrag{-2n+1}{$-2n+1$} 
\psfrag{-n+1}{$\kern-5pt-n+1$} 

{\centering \includegraphics{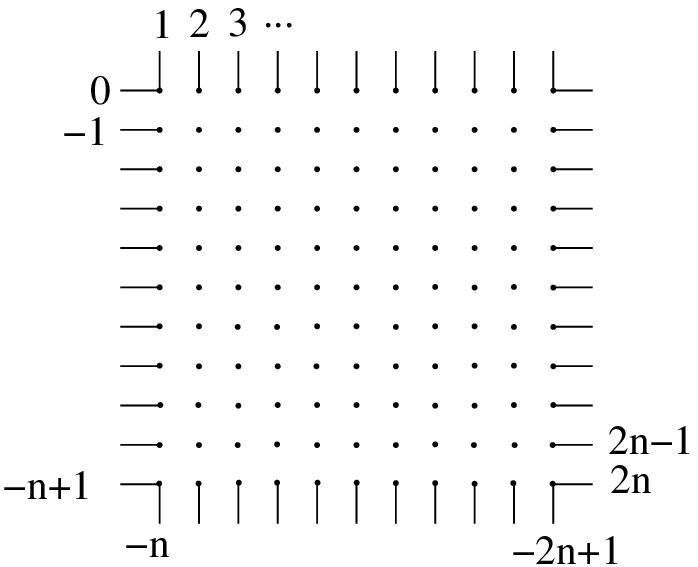} \par}

\caption{The labelling of the external links\label{fpl1}} \end{figure}

The reader
should recall from the introduction that any FPL configuration defines
a matching on the external links occupied by the polygons, by matching
those which are on the same polygon. We call this matching the \emph{matching
associated to the FPL configuration}. When we think of the matching
as being fixed, and when we consider all FPL configurations having
this matching as associated matching, we shall also speak of these
FPL configurations as the ``FPL configurations corresponding to this
fixed matching.''

We label the $4n$ external links around $Q_{n}$ 
by $\{-2n+1,-2n+2,\ldots,2n-1,2n\}$
clockwise starting from the right-most link on the bottom side of
the square, see Figure~\ref{fpl1}. If \( \al \) is an external link
of the square, we denote its label by \( L(\al) \). Throughout this
paper, all the FPL configurations that are considered are configurations
which correspond to matchings of either the even labelled external
links or the odd labelled external links.

\subsection{Wieland's rotational invariance}
Let \( X \) be a non-crossing
matching of the set of even (odd) labelled external
links. Let \( \tilde{X} \) be the ``rotated'' matching of
the odd (even) external links defined by the property that the links
labelled $i$ and $j$ in $X$ are matched if and only if the links
labelled $i+1$ and $j+1$ are matched in $\tilde{X}$, where we identify
$2n+1$ and $-2n+1$. Let \( FPL(X) \)
denote the set of FPL configurations corresponding to the matching
\( X \). Wieland \cite{wie} proved the following surprising result. 

\begin{thm}
[\sc Wieland] \label{thm:wie} For any matching \( X \) of the even
(odd) labelled external links, we have \[ |FPL(X)|=|FPL(\tilde{X})|.\]
\end{thm}
In other terms, the number of FPL configurations corresponding to
a given matching is invariant under rotation of the ``positioning''
of the matching around the square. As we mentioned already in the
introduction, this being the case, we can represent
matchings in terms of chord diagrams of $2n$ points placed around
a circle (see Figure~\ref{fpl25}).

\subsection{Partitions and Ferrers diagrams}\label{sec:hook}
Next we explain our notation concerning partitions and Ferrers
diagrams (see e.g.\ \cite[Ch.~7]{StanBI}).
A \emph{partition} is a vector $\lambda=(\lambda_1,\lambda_2,\dots,\lambda_\ell)$
of positive integers such that $\lambda_1\ge\lambda_2\ge\dots\ge\lambda_\ell$.
For convenience,
we shall sometimes use exponential notation. For example, the partition
$(3,3,3,2,1,1)$ will also be denoted as $(3^3,2,1^2)$. To each partition
$\lambda$, one associates its {\it Ferrers diagram}, which is the
left-justified arrangement of cells with $\lambda_i$ cells in the
$i$-th row, $i=1,2,\dots,\ell$. See Figure~\ref{beau13} for the
Ferrers diagram of the partition $(7,5,2,2,1,1)$. (At this point,
the labels should be disregarded.) We will usually identify a Ferrers
diagram with the corresponding partition; for example we will say
{}``the Ferrers diagram $(\lambda_{1},\ldots,\lambda_{\ell})$''
to mean {}``the Ferrers diagram corresponding to the partition $(\lambda_{1},\ldots,\lambda_{\ell})$''.
The size $|\la|$ of a Ferrers diagram $\la$ is given by the total number
of cells of $\la$. The partition {\it conjugate to} $\lambda$ is the
partition $\lambda'=(\lambda_1',\lambda_2',\dots,\lambda'_{\lambda_1})$,
where $\lambda_j'$ is the length of the $j$-th column of the Ferrers
diagram of $\lambda$. 

Given a Ferrers diagram $\la$,
we write $(i,j)$ for the cell in the $i$-th row and the $j$-th
column of $\la$. We use the notation \( u=(i,j)\in \la \) to express
the fact that $u$ is a cell of $\la$. Given a cell $u$, we denote
by \( c(u):=j-i \) the \emph{content} of \( u \) and by \( h(u):=\lambda _{i}+\lambda'_{j}-i-j+1 \) the \emph{hook length} of \( u \),
where $\lambda$ is the partition associated to $\la$.

It is well-known (see e.g.\ \cite[p.~50]{FuHaAA}), that 
the dimension of the irreducible representation of the symmetric group
$S_{|\la|}$ indexed by a partition (or, equivalently, by a 
Ferrers diagram) $\la$, 
which we denote by $\dim(\la)$, is given by the hook-length formula due
to Frame, Robinson and Thrall \cite{FrRTAA},
\begin{equation} \label{eq:hook} 
\dim(\la)=
\frac {n!} {\prod _{u\in \la} ^{}h_u}.
\end{equation} 

\subsection{How to associate a Ferrers diagram to a
matching}\label{sec:Ferrers}
\begin{figure}[htb]
\begin{center}\includegraphics{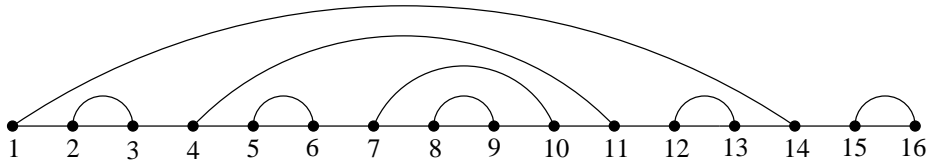}\end{center}

\caption{\label{beau6}A planar matching}
\psfrag{1}{$1$} \psfrag{2}{$2$}\psfrag{3}{$3$} \psfrag{4}{$4$}\psfrag{5}{$5$}\psfrag{6}{$6$}\psfrag{7}{$7$}\psfrag{8}{$8$}\psfrag{9}{$9$}\psfrag{10}{$10$}\psfrag{11}{$11$}\psfrag{12}{$12$}\psfrag{15}{$15$}\psfrag{16}{$16$}
\end{figure}

Let $X$ be a non-crossing matching on the set $\{1,2,\ldots, 2d\}$,
that is, an involution of this set with no fixed points which can be
represented by non-crossing arches in the upper half-plane 
(see Figure~\ref{beau6}
for an example of a non-crossing matching of the set $\{1,2,\ldots,16\}$).
Such a non-crossing matching can be translated into a 0-1-sequence
$v(X)=$ $v_{1}v_{2}\ldots v_{2d}$ of length $2d$
by letting $v_{i}=0$ if $X(i)>i$, and $v_i=1$ if $X(i)<i$.
For example, if $X$ is the matching appearing in Figure~\ref{beau6},
then $v=0010010011101101$.

\begin{figure}[htb]
\psfrag{1}{$1$} \psfrag{0}{$0$} 
\begin{center}\resizebox*{4cm}{!}{\includegraphics{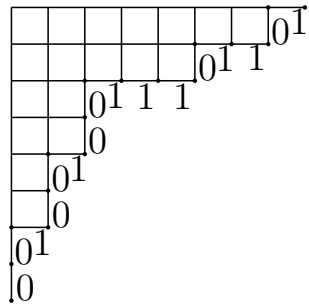}}\end{center}

\caption{\label{beau13}A Ferrers diagram and its $d$-code}
\end{figure}

On the other hand, any 0-1-sequence can be translated into a Ferrers
diagram by reading the 0-1-sequence from left to right and
interpreting a 0 as a unit up-step and a 1 as unit right-step.
{}From the starting point of the obtained path we draw a vertical
segment up-wards, and from the end point a horizontal segment
left-wards. By definition, the region enclosed by the path, the
vertical and the horizontal segment is the {\it Ferrers diagram
associated to the given matching}. See Figure~\ref{beau13} for the
Ferrers diagram associated to the matching in Figure~\ref{beau6}. 
(In the figure, for the sake of clarity, we have labelled the up-steps
of the path by 0 and its right-steps by 1.)
In the sequel, we shall denote the Ferrers diagram associated to $X$
by $\la(X)$. 

Conversely, given a Ferrers diagram $\la$, there are several 0-1-sequences
which produce $\la$ by the above described procedure. 
Namely, by moving along the
lower/right boundary of $\la$ from lower-left to top-right, and
recording a 0 for every up-step and a 1 for every right-step, we
obtain one such 0-1-sequence. Prepending an arbitrary number of 0's and
appending an arbitrary number of 1's we obtain all the other sequences
which give rise to $\la$ by the above procedure. Out of those, we shall
make particular use of the so-called {\it $d$-code} of $\la$ 
(see \cite[Ex.~7.59]{StanBI}). Here, $d$ is a positive integer such
that $\la$ is contained in the Ferrers
diagram $(d^{d})$. We embed $\la$ in $(d^d)$ so that the diagram $\la$
is located in the top-left corner of the square $(d^d)$. 
We delete the lower side and the right side of the square $(d^d)$.
(See Figure~\ref{beau13} for an example
where $d=8$ and $\la=(7,5,2,2,1,1)$.) 
Now, starting from the lower/left corner of the square, we move, as
before, along the
lower/right boundary of the figure from lower-left to top-right,
recording a 0 for every up-step and a 1 for every right-step.
By definition, the obtained 0-1-sequence is the $d$-code of $\la$. 
Clearly, the $d$-code has exactly $d$ occurrences of $0$
and $d$ occurrences of $1$. 
For example, the $8$-code of the Ferrers diagram $(7,5,2,2,1,1)$
is $0010010011101101$.

\subsection{\label{secHook}An enumeration result for rhombus tilings}
In the proof of Theorem~\ref{main}, 
we shall need a general result on the enumeration of
rhombus tilings of certain subregions
of the regular triangular lattice in the plane,
which are indexed by Ferrers diagrams. 
(Here, and in the sequel,
by a rhombus tiling we mean a tiling by rhombi of unit side lengths
and angles of $60^\circ$ and $120^\circ$.) 
This result appeared in an equivalent form in
\cite[Theorem~2.6]{caskra}. As is shown there, it follows from 
Stanley's hook-content formula \cite[Theorem~15.3]{StanAA}, via
the standard bijection between rhombus tilings and non-intersecting lattice
paths, followed by the standard 
bijection between non-intersecting lattice paths
and semistandard tableaux.

\begin{figure}[htb]
 \psfrag{d}{$d$} 
\psfrag{h}{$h$} 
\psfrag{d+h}{$d+h$} 
\psfrag{0}{$0$} 
\psfrag{1}{$1$} 
\psfrag{d-code of X}{$d$-code of $\la$}
\includegraphics{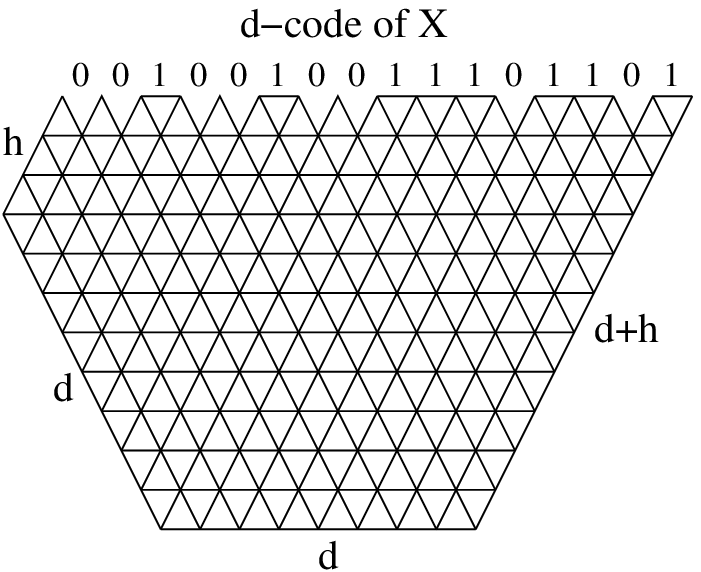}

\caption{\label{beau14}The region $R(\la,d,h)$}
\end{figure}

Let $\la$ be a Ferrers diagram contained in the square $(d^{d})$,
and let $h$ be a non-negative integer $h$.
We define the region \( R(\la,d,h) \)
to be a pentagon with some notches along the top side. More precisely
(see Figure~\ref{beau14} where the region  \( R(\la,8,3) \) is shown,
with $\la$ the Ferrers diagram $(7,5,2,2,1,1)$ from Figure~\ref{beau13}),
the region $R(\la,d,h)$ is the pentagon with base side and bottom-left
side equal to $d$, top-left side $h$, a top side of length $2d$
with notches which will be explained in just a moment, and right side
equal to $d+h$. To determine the notches along the top side, we read 
the $d$-code of $\la$, and we
put a notch whenever we read a $0$, while we leave a horizontal
piece whenever we read a $1$.

We can now state the announced enumeration result for rhombus tilings
of the regions $R(\la,d,h)$. 

\begin{thm}
\label{thm:rhla}Given a Ferrers diagram $\la$ contained in the square
$(d^{d})$ and a positive integer $h$, the number of
rhombus tilings of the region $R(\la,d,h)$ is given by
$SSYT(\la,d+h)$,
where
\begin{equation} \label{eq:hook-content} 
SSYT(\la,N)=
\prod _{u\in \la} ^{}\frac {c(u)+N} {h(u)},
\end{equation}
with $c(u)$ and $h(u)$ the content and the hook length of $u$,
respectively, as defined in Section~{\em\ref{sec:hook}}.
\end{thm}
\begin{figure}[htb]
\psfrag{h1}{}
\psfrag{h2}{}
\psfrag{v1}{}
\psfrag{v2}{}
\psfrag{hk}{}
\psfrag{vk}{}
\psfrag{r}{}
\psfrag{l1}{}
\psfrag{l1'}{}
\psfrag{r+}{}

{\centering \includegraphics{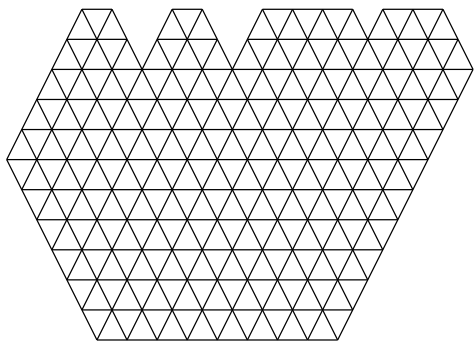} \par}

\caption{The reduced region\label{beau32}}
\end{figure}

\begin{remark}
The choice of the notation $SSYT(\la,N)$ comes from the fact that
the number in \eqref{eq:hook-content} counts at the same time the
number of semistandard tableaux of shape $\la$ with entries at most
$N$ (cf.\ \cite[Theorem~15.3]{StanAA}). Indeed, implicitly in the proof of the
above theorem which we give below is a bijection between the rhombus
tilings in the statement of the theorem and these semistandard tableaux.
\end{remark}
\begin{proof}[Proof of Theorem~\ref{thm:rhla}]
It should be observed that, due to the nature of the
region\break $R(\la,d,h)$, there are several ``forced'' subregions, that is,
subregions where the tiling is completely determined. For example, 
the right-most layer in Figure~\ref{beau14} must necessarily be
completely filled with right-oriented rhombi, while
the first two upper-left layers must be filled with horizontally
symmetric rhombi. If we remove all the ``forced'' rhombi,
then a smaller region remains. See Figure~\ref{beau32} for the result
of this reduction applied to the region in Figure~\ref{beau14}.
To the obtained region we may apply Theorems~2.6 and 2.5 
from \cite{caskra}. As a result, we obtain the
desired formula.
\end{proof}

\subsection{The Lindstr\"om--Gessel--Viennot formula}
It is well-known that rhombus tilings are (usually) in bijection with
families of non-intersecting lattice paths. We shall make use of this
bijection in Section~\ref{secGen}, together with
the main result on the enumeration of non-intersecting lattice
paths, which is a determinantal formula due to Lindstr\"om
\cite{LindAA}.
In the combinatorial
literature, it is most often attributed to Gessel and Viennot
\cite{GeViAA,gv}, but it can actually be traced back to Karlin and
McGregor \cite{KaMGAB,KaMGAC}.

Let us briefly recall that formula,
or, more precisely, a simplified version tailored for our purposes.
We consider lattice paths in the planar integer lattice $\mathbb Z^2$ 
consisting of unit horizontal and vertical steps in the positive
direction. Given two points $A$ and $E$ in $\mathbb Z^2$, 
we write $\P (A \to E)$ for the number of
paths starting at $A$ and ending at $E$.
We say that a family of paths is {\it non-intersecting} if no two
paths in the family have a point in common.

We can now state the announced main result on non-intersecting lattice paths
(see \cite[Lemma~1]{LindAA} or \cite[Cor.~2]{gv}). 
\begin{lem} \label{gv}
Let $A_1,A_2,\dots, A_n, E_1, E_2,\dots , E_n$ be points of 
the planar integer lattice $\mathbb Z^2$, such that for all $i<j$ the point
$A_i$ is {\rm(}weakly{\rm)} south-east of the point $A_j$, and the
point $E_i$ is {\rm(}weakly{\rm)} south-east of the point $E_j$.
Then the number of families $(P_1,P_2,\dots,P_n)$ of non-intersecting
lattice paths, $P_i$ running from $A_i$ to $E_i$, $i=1,2,\dots,n$, is
given by
\begin{equation} \label{eq:gv}
\det{\left(\P (A_i \to E_j)\right)}_{1\le i,j \le n}.
\end{equation}
\quad \quad \qed
\end{lem}

\section{\label{secfix}Fixed edges}

In this section, we perform the first step in order to prove
Conjecture~\ref{conj:zub1}. Let $X$ be a given non-crossing 
matching with $d$ arches.
As in the statement of the conjecture, let $X\cup m$ be the matching
obtained by adding $m$ nested arches to $X$.
Thanks to Theorem~\ref{thm:wie} of
Wieland, we may place $X\cup m$ in an arbitrary way around $Q_n=Q_{d+m}$, 
without changing the number of corresponding FPL configurations.
We place $X\cup m$ so that, using Lemma~\ref{fixedg} below, the FPL
configurations corresponding to the matching will have as many forced
edges as possible. 

\begin{figure}
 \psfrag{T}{$X$}
 \psfrag{M}{$M$}
\begin{center}\includegraphics{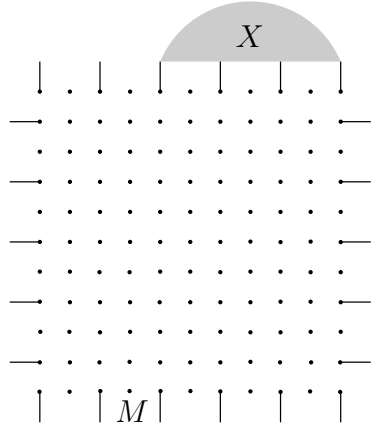}\end{center}
\caption{\label{beau10}Placing the matching around $Q_{n}$}
\end{figure}

\begin{figure}
 \psfrag{M}{\Large$M$}
 \psfrag{T}{\huge$X$}
 \psfrag{A}{\Large$A$}
 \psfrag{B}{\Large$B$}
 \psfrag{C}{\Large$C$}
 \psfrag{D}{\Large$D$}
 \psfrag{E}{\Large$E$}
 \psfrag{J}{\Large$J$}
 \psfrag{K}{\Large$K$}
 \psfrag{L}{\Large$L$}
 \psfrag{F}{\Large$F$}
 \psfrag{G}{\Large$G$}
\begin{center}\includegraphics[%
  width=0.95\textwidth]{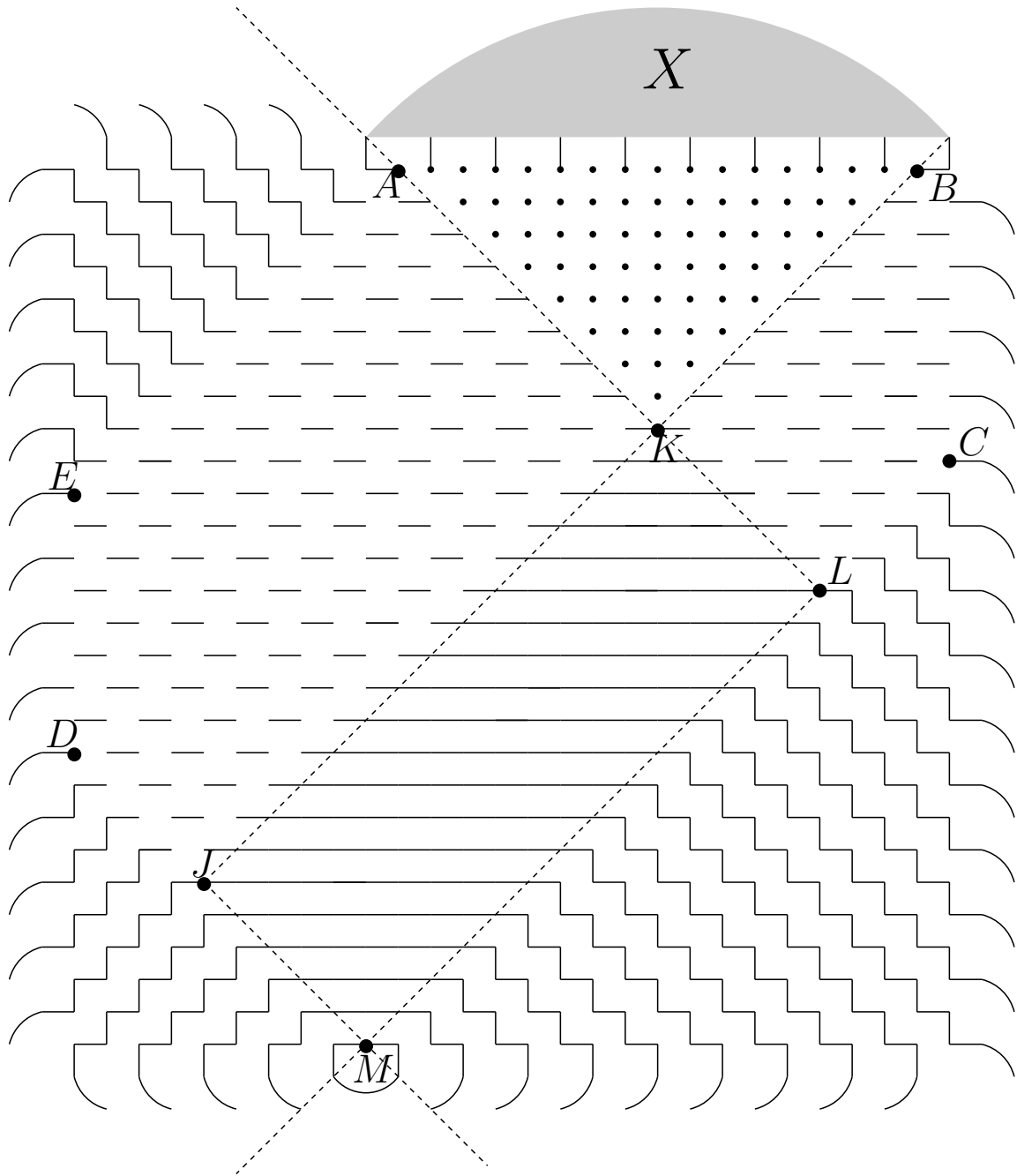}\end{center}

\caption{\label{beau12}The set of fixed edges}
\end{figure}

To be precise,
we place $X\cup m$ so that the arches corresponding to $X$ appear on the
very right of the upper side of $Q_n$.
That is, we place these arches on the external links labelled 
$ n-4d+2, n-4d+4,\ldots, n-2, n$. 
Equivalently, we choose to place the centre $M$ of the $m$
nested arches on the external link labelled by $-n-2d+1$.
See Figure~\ref{beau10} for a
schematic picture, and Figure~\ref{beau12} for a more elaborate one (in
which the added edges in the interior of the square grid $Q_n$ should
be ignored at this point).
In order to guarantee that $X$ has
place along the upper side of the square grid, we must assume that
$m\ge 3d$.

The following lemma helps to identify edges which are occupied by {\it
each} FPL configuration corresponding to a given matching.
It is a consequence of an iterated use of a result
of de Gier (see \cite[Lemma~8]{degier} or 
\cite[Lemmas~2.2 and 2.3]{caskra}). In the sequel,
when we speak of ``fixed edges'' we always mean edges that
have to be occupied by \emph{any} FPL configuration under consideration.

\begin{lem}
\label{fixedg}Let \( \al =\al _{1},\al _{2},\ldots ,\al _{k}=\be \) be a sequence
of external links, where $L(\al _i)=a+2i$ \emph{mod} $4n$, for some
fixed $a$, that is, the external links $\al _1,\al _2,\dots,\al _k$ comprise
every second external link along the stretch between $\al $ and $\be$
along the border of $Q_n$ \emph{(}clockwise\emph{)}. Furthermore,
we suppose that one of the following conditions holds: 
\begin{enumerate}
\item \( \al  \) and \( \be \) are both on the top side of \( Q_{n} \), that
is, \( 1\leq L(\al )<L(\be)\leq n \);
\item \( \al  \) is on the top side and \( \be \) is on the right side of \( Q_{n} \),
that is, \( 1\leq L(\al )\leq n<L(\be) \) and \( n-L(\al )>L(\be)-(n+1) \);
\item \( \al  \) is on the left side and \( \be \) is on the right side of
\( Q_{n} \), that is, \( n<L(\be)\leq 2n \) and \( -n<L(\al )\leq 0 \).
\end{enumerate}
\begin{figure}

\psfrag{O}{\Large$O$}

\psfrag{A}{\Large$\al $}

\psfrag{B}{\Large$\be$}

\psfrag{A'}{\Large$A'$} 

\psfrag{B'}{\Large$B'$}

{\centering {\includegraphics[%
  width=0.85\textwidth]{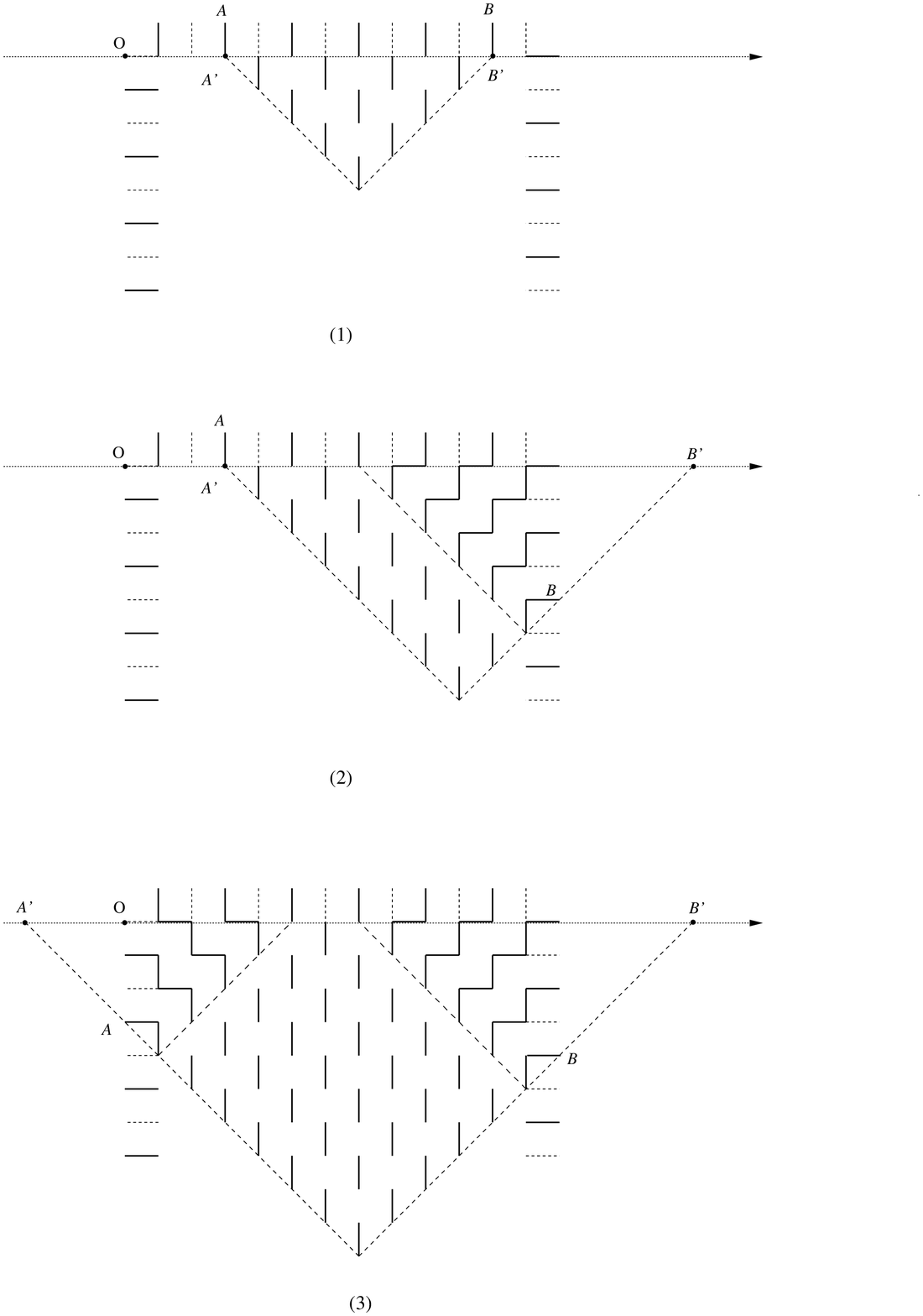}} \par}

\caption{The possible regions of fixed edges determined by a sequence of external links belonging to distinct loops\label{seq}} 

\end{figure}

For the FPL configurations for which the external links 
\( \al _{1},\al _{2},\ldots ,\al _{k} \)
belong to different loops, the region of fixed edges is
\emph{(}essentially\emph{)}
triangular \emph{(}see Figure~\emph{\ref{seq}} for illustrations
of the region and the fixed edges in its interior; the ``essentially''
refers to the fact that in Cases~\emph{(2)} and \emph{(3)} parts
of the triangle are cut off\emph{)}. More precisely, if one places
the origin \( O \) of the coordinate system one unit to the left
of the top-left corner of \( Q_{n} \), the coordinates of the triangle
are given in the following way: let \( A' \) and \( B' \) be the
points on the $x$-axis with \( x \)-coordinates \( L(\al ) \) and
\( L(\be) \), respectively, then the region of fixed edges is given
by the intersection of the square \( Q_{n} \) and the \emph{(}rectangular
isosceles\emph{)} triangle having the segment \( A'B' \) as basis.

In Cases~\emph{(2)} and \emph{(3)}, the configurations are completely
fixed as ``zig-zag'' paths in the corner regions of \( Q_{n} \)
where a part of the triangle was cut off \emph{(}see again Figure~\emph{\ref{seq})}.
More precisely, in Case~\emph{(2)}, this region is the reflexion
of the corresponding cut off part of the triangle in the right side
of \( Q_{n} \), and in Case~\emph{(3)} it is that region and also
the reflexion of the corresponding cut off part on the left in the
left side of \( Q_{n} \).
\end{lem}

We use this lemma to determine the set of fixed edges of
the FPL configurations corresponding to the matching $X\cup m$ 
in Conjecture~\ref{conj:zub1}.
For convenience (the reader should consult Figure~\ref{beau12} while
reading the following definitions), we let $A,B,C,D,E$ be the border vertices 
of the external links labelled $n-4d+3$, $n-1$, $n+2d$, $-n+2d$,
$-n+4d-2$, respectively, we let
$J$ be the intersection point of the line connecting $D$ and $M$ and
the line emanating diagonally from $B$, we let $K$ be the intersection
point of the latter line emanating from $B$ and the line emanating
diagonally (to the right) from $A$, and we let $L$ be the intersection
point of the latter line emanating from $A$ and the line connecting
$C$ and $M$.
We state the result of the application of Lemma~\ref{fixedg} to our case
in form of the following lemma.

\begin{lem}
The region of fixed edges of the FPL configurations corresponding
to the matching in Conjecture~{\em\ref{conj:zub1}} contains all the
edges indicated in Figure~\emph{\ref{beau12}}, that is:

\begin{enumerate} 
\item all the horizontal edges in the rectangular region $JKLM$,
\item every other horizontal edge in the pentagonal region $AKJDE$ as
indicated in the figure,
\item every other horizontal edge in the region $BCLK$ as
indicated in the figure,
\item the zig-zag lines in the corner regions above the line $EA$,
respectively below the lines $DM$ and $MC$, as indicated in the figure.
\end{enumerate}

\end{lem}
\begin{proof}
This result follows by applying Lemma~\ref{fixedg} to
all the external links corresponding to the $m$ nested arches
which are on the left of the centre $M$ plus the {}``first''
external link of the matching $X$, on the one hand, 
and to all the external links corresponding
to the $m$ nested arches which are on the right of the centre $M$
plus the {}``last'' external link of the matching 
$X$, on the other hand. More precisely, we apply 
Lemma~\ref{fixedg} to the sets 
$$\{ \al\text{ an external link with }-n-2d+2\le L(\al)\leq 
n-4d+2\textrm{ or }L(\al)\geq3n-2d\}$$
and 
$$\{ \al\text{ an external link with either }L(\al)\leq-n-2d\text{ or
}L(\al)\ge n\}.$$
The triangles forming the respective regions of fixed edges
are drawn by dashed lines in Figure~\ref{beau12}. 
Note that the two triangles of fixed
edges overlap in the rectangular region $JKLM$,
where the fixed edges form
parallel horizontal lines.
\end{proof}

\section{\label{secCon} Proof of Conjecture~\ref{conj:zub1} for $m$
large enough}

Let $X$ be a non-crossing 
matching consisting of $d$ arches, and let $m$ be a positive
integer such that $m\geq3d$. As in Conjecture~\ref{conj:zub1}, we
denote the matching arising by adding $m$ nested arches to $X$ by $X\cup m$.
The goal of this section is to give
an explicit formula for the number $A_{X}(m)$ of FPL configurations
corresponding to the matching $X\cup m$
(cf.\ Figure~\ref{beau8}), which implies Conjecture~\ref{conj:zub1}
for $m\ge 3d$, see
Theorem~\ref{main}. 

\begin{figure}
\psfrag{T}{
\huge$X$} 
\psfrag{M}{\Large$M$} 
 \psfrag{A}{\Large$A$}
 \psfrag{B}{\Large$B$}
 \psfrag{C}{\Large$C$}
 \psfrag{D}{\Large$D$}
 \psfrag{E}{\Large$E$}
 \psfrag{J}{\Large$J$}
 \psfrag{K}{\Large$K$}
 \psfrag{L}{\Large$L$}
 \psfrag{F}{\Large$F$}
 \psfrag{G}{\Large$G$}
\psfrag{xi}{\huge$\xi_1$}
\includegraphics[width=1.0\textwidth]{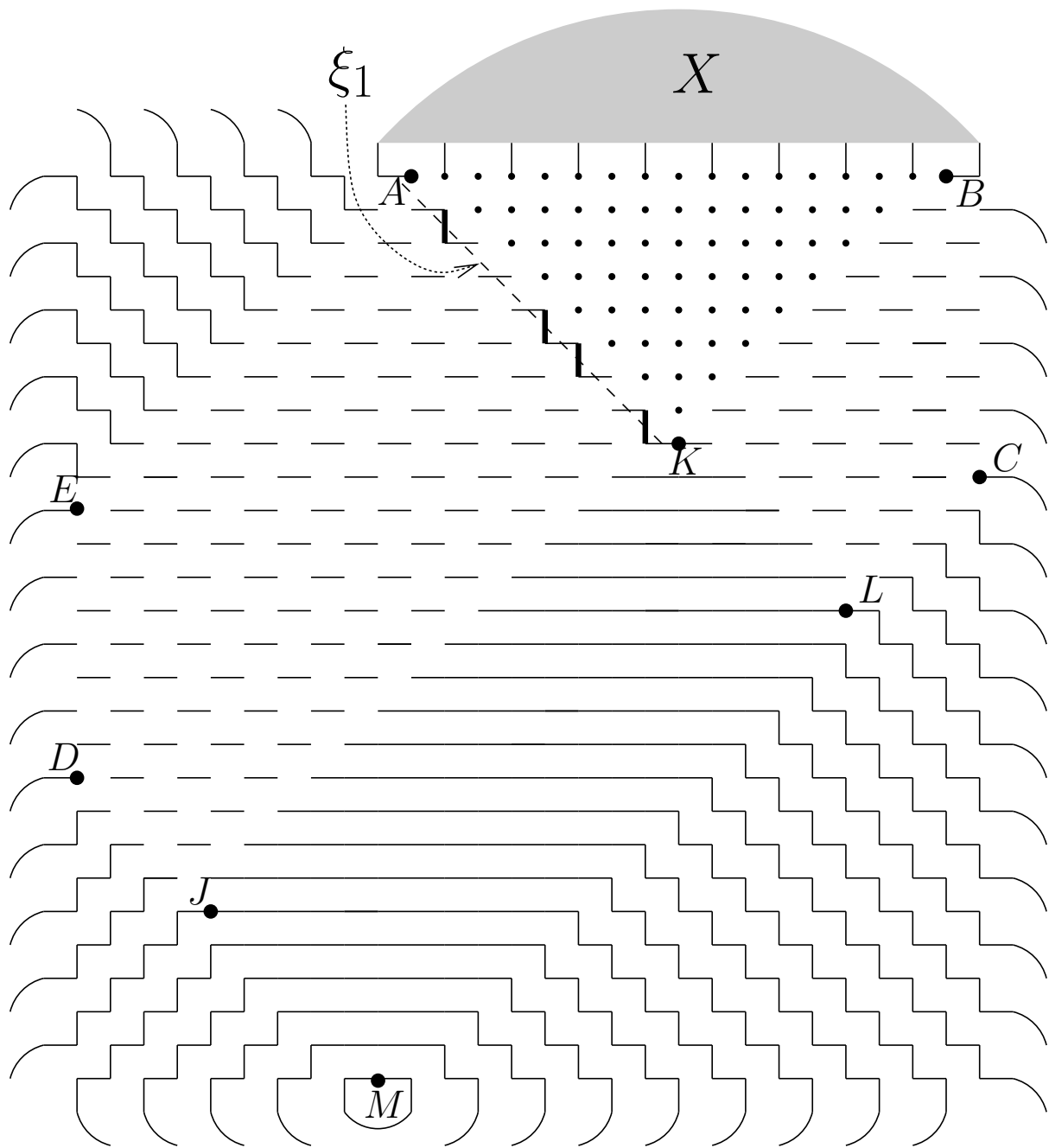}

\caption{\label{beau15}Splitting the problem}
\end{figure}

Recall
from Section~\ref{secfix} how we place this matching around the grid
$Q_n=Q_{d+m}$, see Figure~\ref{beau12}, where we have chosen $d=5$, $m=23$,
and, hence, $n=d+m=28$. The reader should furthermore
recall the placement of the points $A,B,C,D,E,J,K,L,M$.
Let $\xi_1$ be the segment which
connects the point which is half a unit to the left of $A$ and the point
which is half a unit left of $K$, see Figure~\ref{beau15}.

There are exactly $2d-2$ possible vertical edges that cross
the segment $\xi_1$. We denote them by $e_{1},e_{2},\ldots,e_{2d-2}$,
starting from the top one and proceeding south-east. 
We claim that among those there are exactly $d-1$ which are occupied
by an FPL configuration. To see this, we first note that in the
rectangular region $JKLM$ there are $m-d+1$ parallel lines strictly below
$K$. Each of them must be part of one of the $m+1$ 
loops starting on the external links labelled
$\{-n-2d+2,-n-2d+4,\ldots,n-4d,n-4d+2\}$ (these are the external
links ``between'' $M$ and $A$, in clockwise direction).
Hence there are exactly $d$ loops that cross the segment $\xi_1$, which
implies that any FPL configuration occupies exactly $d-1$
vertical edges out of $\{e_{1},e_{2},\ldots,e_{2d-2}\}$, as we claimed. 
We encode a choice of $d-1$ edges from
$\{e_{1},e_{2},\ldots,e_{2d-2}\}$ by a subset $\mathcal E$ from 
$\{1,2,\dots,2d-2\}$, by making the obvious identification that the
choice of $e_{i_1},e_{i_2},\dots,e_{i_{d-1}}$ is encoded by $\mathcal
E=\{i_1,i_2,\dots,i_{d-1}\}$.

On the other hand, any such choice
is equivalent to the choice of a Ferrers diagram contained in the
square Ferrers diagram $((d-1)^{d-1})$ by the following construction. Let
$\mathcal E\subset\{ {1},{2},\ldots,{2d-2}\}$ be of cardinality $d-1$.
Let $c_{\mathcal E}:=c_{1}c_{2}\ldots c_{2d-2}$ be the binary string defined
by $c_{i}=0$ if and only if ${i}\in \mathcal E$. 
The string $c_{\mathcal E}$ obtained
in this way determines a Ferrers diagram, as we described in
Section~\ref{sec:hook}. We denote this Ferrers diagram by $\la(\mathcal E)$.
In the example in Figure~\ref{beau15}, 
we have $\mathcal E=\{ {2},{5},{6},{8}\}$
and, hence, $\la(\mathcal E)=(4,3,3,1)$.

It is obvious from the picture, that, once we have chosen the vertical
edges along $\xi_1$  
belonging to an FPL configuration with associated matching $X\cup m$
(that is, the vertical edges out of
$\{e_1,e_2,\dots,e_{2d-2}\}$ which are occupied by the FPL
configuration), the configuration can be completed
separately in the region to the ``left'' of $\xi_1$ (that is, in the
region $AKJDE$) and to the ``right'' of $\xi_1$ (that is, in the region
$ABCLK$). In particular, it is not difficult to see 
that that the number of FPL configurations
with associated matching $X\cup m$ 
which, out of $\{e_1,e_2,\dots,e_{2d-2}\}$, occupy a {\it fixed\/} subset of
vertical edges, is equal to the number of FPL configurations in the
region $AKJDE$  
times the number of FPL configurations in the region
$ABCLK$ which respect the matching $X$.

\begin{figure}
\psfrag{xi1}{\Large$\xi_1$} 
\psfrag{xi2}{\Large$\xi_2$}
\psfrag{X}{\Large$X$}
\psfrag{X'}{\Large$\mathcal E$}
\includegraphics{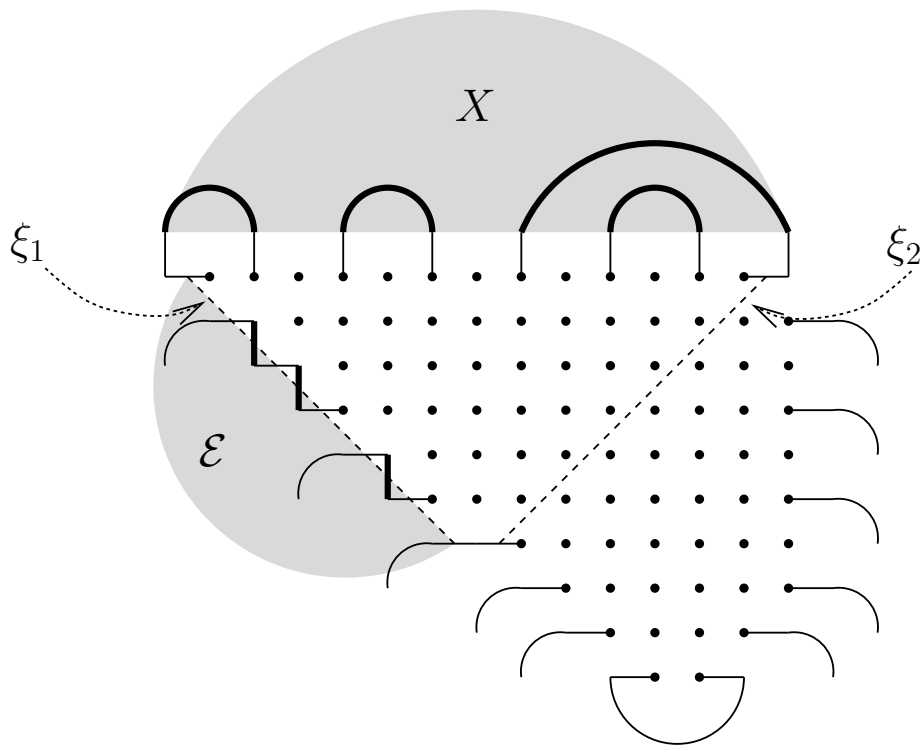}
\caption{\label{beau16}The numbers $a_{X}(\mathcal E)$}
\end{figure}
\begin{figure}
\includegraphics[%
  width=1.0\textwidth]{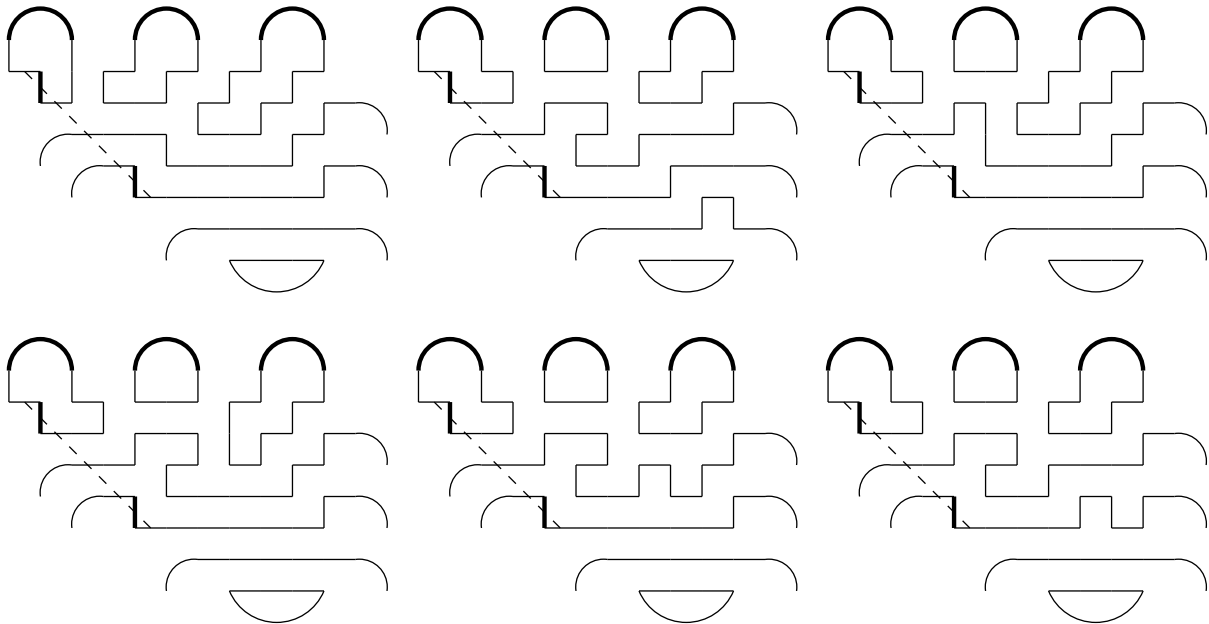}
\caption{\label{beau18}The six configurations corresponding to 
$X=\{1\leftrightarrow 2,\,3\leftrightarrow 4,\,5\leftrightarrow 6\}$ and
$\mathcal E=\{1,4\}$}
\end{figure}

Clearly, the region $ABCLK$ to the right of $\xi_1$ does not depend on
$m$. We denote the number of FPL configurations of that region which
respect the matching $X$ and whose set of edges from
$\{e_1,e_2,\dots,e_{2d-2}\}$ is encoded by $\mathcal E$
by $a_{X}(\mathcal E)$.
For example, if $X$ is the matching
$\{1\leftrightarrow 2,\,3\leftrightarrow 4,\,5\leftrightarrow 6\}$, 
and if $\mathcal E$ is the set $\{1,4\}$,
then we have $a_{X}(\mathcal E)=6$.
The six configurations corresponding to this choice of $X$ and $\mathcal
E$ are
shown in Figure~\ref{beau18}, where the arches corresponding to $X$
and the edges corresponding to $\xi_1$ are marked in bold-face.
We have $\la(X)=(2,1)$ and $\la(\mathcal E)=(1,1)$. In particular,
$\la(\mathcal E)\subseteq \la(X)$. The next lemma shows that this is
not an accident. 

\begin{lem}
\label{imposs}
Let $X$ be a non-crossing matching with $d$ arches and let $\mathcal E$
be a subset of $\{1,2,\dots,{2d-2}\}$ consisting of $d-1$ elements.

\begin{enumerate}
\item If
$\la(\mathcal E)\not\subseteq \la(X)$, then $a_{X}(\mathcal E)=0$.
\item
If $\la(\mathcal E)=\la(X)$, then $a_{X}(\mathcal E)=1$.
\end{enumerate}
\end{lem}

\begin{proof}This follows from Corollary~\ref{cor}.(1) and (3).
\end{proof}

Equipped with this lemma, we are now able to prove the first main
result of this paper.

\begin{thm}
\label{main}Let $X$ be a non-crossing matching with $d$ arches, 
and let $m\geq3d$.
Then
\begin{equation} \label{eq:main} 
A_{X}(m)=SSYT(\la(X),m-2d+1)+
\sum_{\mathcal E:\la(\mathcal E)\varsubsetneq \la(X)}
a_{X}(\mathcal E)\cdot SSYT(\la(\mathcal E),m-2d+1).
\end{equation}
\end{thm}
\begin{proof}
Let us fix $d-1$ edges from
$\{e_1,e_2,\dots,e_{2d-2}\}$, encoded by the set $\mathcal E$. 
In view of Lemma~\ref{imposs}, 
it suffices to show that the number
of configurations on the left of the segment $\xi_1$
(more precisely,
in the region $AKJDE$; see Figure~\ref{beau15}) which, out of
$\{e_1,e_2,\dots,e_{2d-2}\}$, 
occupy exactly the edges encoded by $\mathcal E$
is equal to $SSYT(\la(\mathcal
E),m-2d+1)$. To do so,
we proceed in a way similar to the proof of the main results in
\cite{caskra}. That is, we translate the problem of enumerating the
latter FPL configurations into a problem of enumerating rhombus
tilings. 

We say that a vertex is {\it free} if it belongs to {\it exactly
one} fixed edge. We draw a triangle around any free vertex in the region
$AKJDE$ in such a way that two free vertices are neighbours
if and only if the corresponding triangles share an edge. In the case
which is illustrated in Figure~\ref{beau15}, this leads to the picture 
in Figure~\ref{beau21}.%

\begin{figure}
\includegraphics{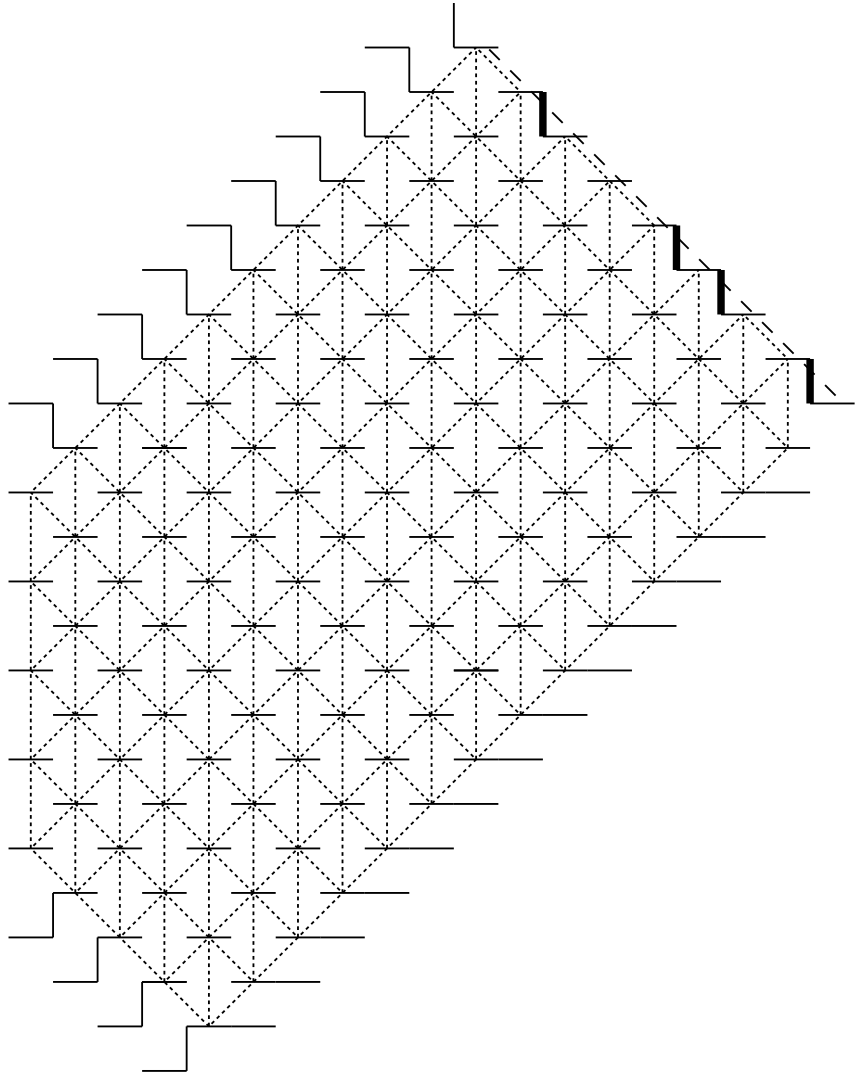}

\caption{\label{beau21}Drawing triangles}
\end{figure}
 Now we make a deformation of the obtained set of triangles in such a way
that all the internal angles become $60^{\circ}$. As a result, we obtain
the region $R(\la(\mathcal E),d-1,m-3d+2)$ defined in
Section~\ref{secHook}, see Figure~\ref{beau22}.%
\begin{figure}
\includegraphics{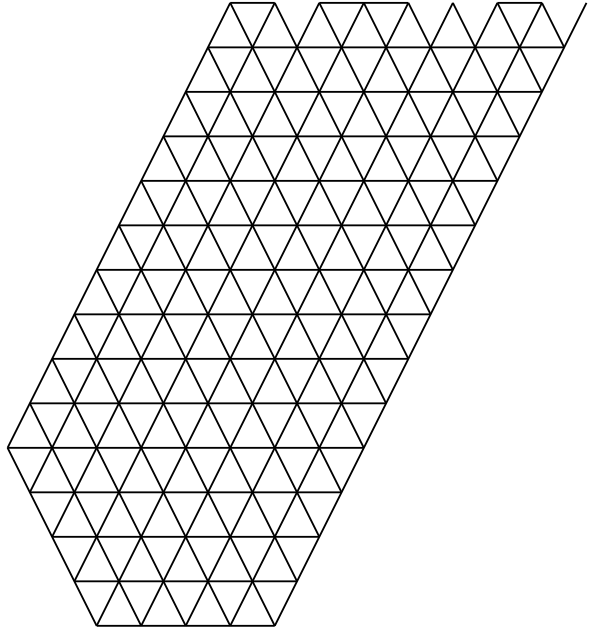}
\caption{\label{beau22}Getting the regions $R(\mathcal E,m+3d-2)$}
\end{figure}
As in \cite{caskra}, it is not difficult to see that the
FPL configurations in the region $AKJDE$
are in bijection with the rhombus tilings of the region 
$R(\la(\mathcal E),d-1,m-3d+2)$. Indeed, to go from a rhombus tiling
to the corresponding FPL configuration, for every rhombus in the
tiling one connects the free vertices which are in the interior of the two
triangles forming the rhombus by an edge. Hence the
result follows from Theorem~\ref{thm:rhla}.
\end{proof}

Zuber's Conjecture~\ref{conj:zub1}, in the case that $m\ge 3d$, 
is now a simple corollary of the above theorem.

\begin{proof}[Proof of Conjecture~\ref{conj:zub1} for $m\ge 3d$]
The polynomiality in $m$ of $A_X(m)$ is obvious from \eqref{eq:main} and
\eqref{eq:hook-content}. The assertion about the integrality of the
coefficients of the ``numerator'' polynomial $P_X(m)$ follows from the
simple fact that the hook product $
\prod _{u\in \la} ^{}h_u$ is a divisor of $\vert\la\vert!$ for any
partition $\la$. Finally,
to see that the leading coefficient of $P_X(m)$ is $\dim(\la(X))$,
one first observes that the leading term in
\eqref{eq:main} appears in the term $SSYT(\la(X),m-2d+1)$. 
The claim follows now by a
combination of \eqref{eq:hook-content} and \eqref{eq:hook}.
\end{proof}

\section{\label{secCon1}Proof of Conjecture~\ref{conj:zub1} for small $m$}

\begin{figure}[h]
\psfrag{T}{\LARGE$X$}\psfrag{S}{$Y$}
\psfrag{xi}{\LARGE$\xi_1$}
\psfrag{x}{\LARGE$\xi_2$}
\resizebox*{15cm}{!}{\includegraphics{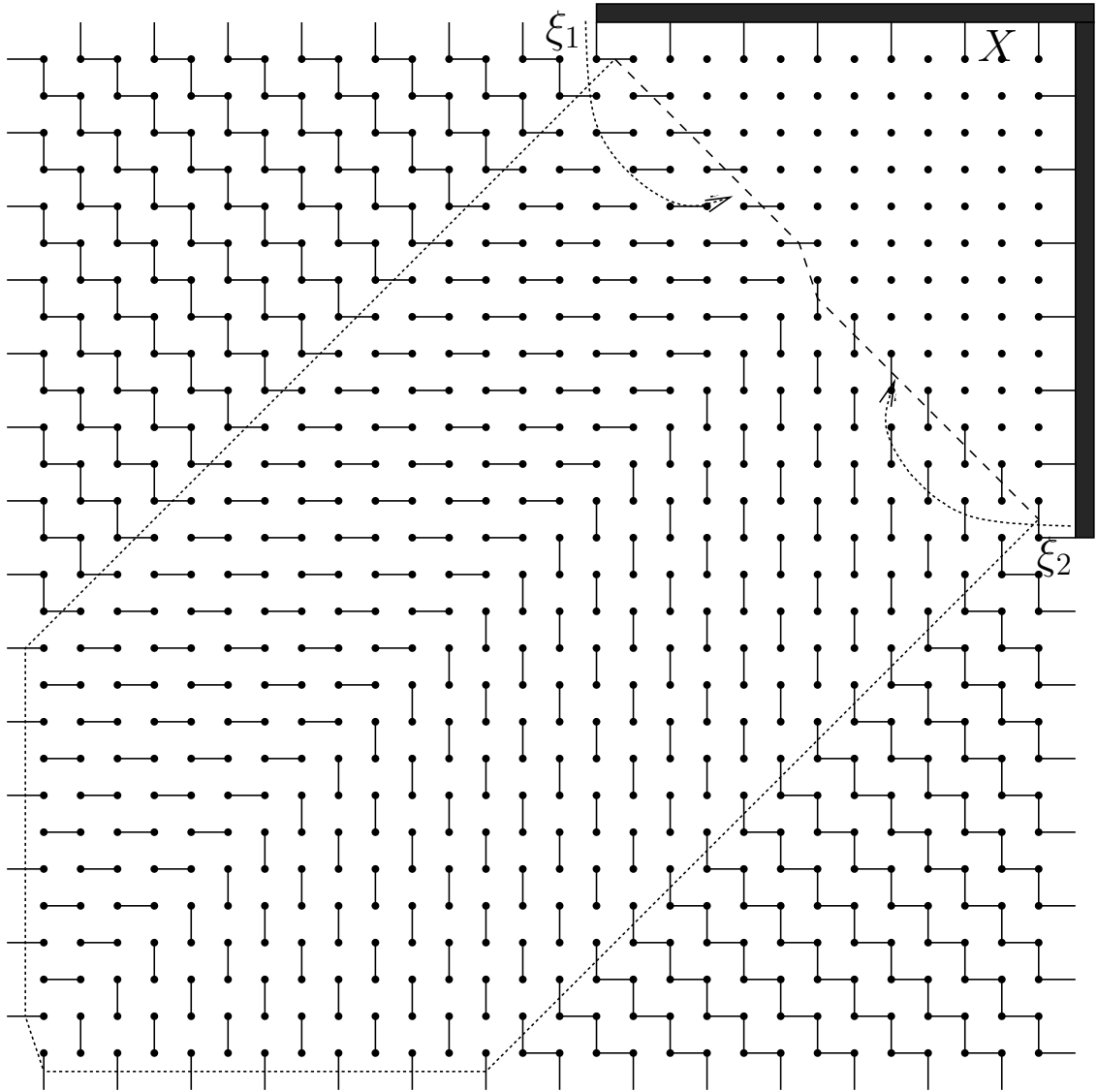}}
\caption{\label{fig:Ecke1}A different placement of the matching around
the grid}
\end{figure}

To prove Conjecture~\ref{conj:zub1} for $m<3d$, we choose a
different placement of the matching $X\cup m$, namely, we place $X$ around the
top-right corner of the square $Q_n$. To be precise, we place $X\cup m$ so
that the arches corresponding to $X$ occupy the external links
labelled $n-2d+2,n-2d+4,\dots,n+2d$, see Figure~\ref{fig:Ecke1}
for an example where $n=28$, $d=7$, and $m=21$. (There, the positioning of the
matching $X$ is indicated by the black hook. The
edges in the interior of the square grid should be ignored at this point.)
In fact, the figure shows an example where $m\ge 2d$, and, strange as
it may seem, this is what we shall assume in the sequel. Only at the
very end, we shall get rid of this assumption.

We now apply again Lemma~\ref{fixedg} to determine the edges which
are occupied by {\it each\/} FPL configuration with associated
matching $X\cup m$. As a result, there are fixed edges along zig-zag
paths in the upper-left and the lower-right corner of the square grid,
while in a pentagonal region located diagonally from lower-left to
upper-right every vertex is on exactly one fixed edge, as
indicated in Figure~\ref{fig:Ecke1}.
(There, the pentagonal region is indicated by the dashed lines.) 
More precisely, the pentagonal region
decomposes into two halves: in the upper-left half every other {\it
horizontal\/} edge is taken by any FPL configuration, whereas in the
lower-right half it is every other {\it vertical\/} edge which is
taken by any FPL configuration. 

\begin{figure}[h]
\psfrag{T}{\LARGE$X$}\psfrag{S}{$Y$}
\psfrag{xi}{\LARGE$\xi_1$}
\psfrag{x}{\LARGE$\xi_2$}
\resizebox*{15cm}{!}{\includegraphics{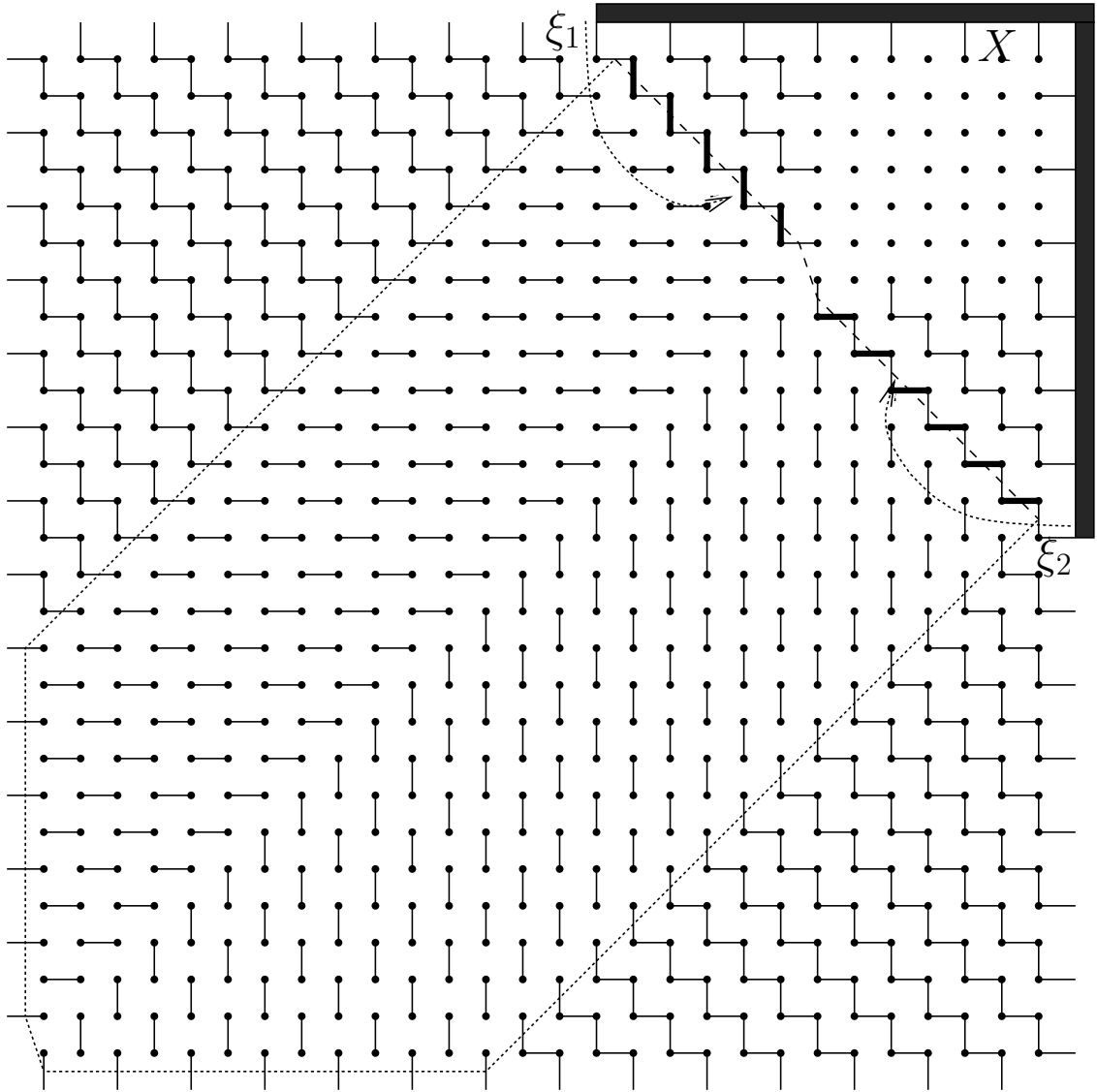}}
\caption{\label{fig:Ecke2}Fixed edges for
$\mathcal E=\{1,\dots,{d-2}\}$ and $\mathcal F=\{1,\dots,{d-1}\}$}
\end{figure}

Now we are argue similarly to the proof of Theorem~\ref{main}.
Along the upper-right border, we mark the segments $\xi_1$ (the
upper-left half; see Figure~\ref{fig:Ecke1}) and $\xi_2$ (the
lower-right half). Let $\{e_1,e_2,\dots,e_{d-2}\}$ be the set of
vertical edges which are crossed by $\xi_1$, and let 
$\{f_1,f_2,\dots,f_{d-1}\}$ be the set of
vertical edges which are crossed by $\xi_2$. In Figure~\ref{fig:Ecke2}
these edges are marked in bold face.
(Compare with Figures~\ref{beau15} and \ref{beau21}.) 
Since any loop which enters the
triangular region in the upper-right corner of the square grid (that
is, the region to the right of the segments $\xi_1$ and $\xi_2$) must
necessarily also leave it again, any FPL configuration with associated
matching $X\cup m$ occupies a subset of
$\{e_1,e_2,\dots,e_{d-2}\}$, encoded by $\mathcal E$ as before, 
and a subset of
$\{f_1,f_2,\dots,f_{d-1}\}$, encoded by $\mathcal F$, such that $\vert\mathcal
E\vert=\vert\mathcal F\vert-1$. Once a choice of $\mathcal E$ and
$\mathcal F$ is made, the number of FPL configurations which cover
exactly the vertical edges encoded by $\mathcal E$ and $\mathcal F$ 
decomposes into the product of
the number of possible configurations in the pentagonal region times
the number of possible configurations in the triangular region in the
upper-right corner of the square grid. Let us denote the former number
by $N(\mathcal E,\mathcal F,m,d)$, and the latter by $c(\mathcal
E,\mathcal F)$. Writing, as before, $A_X(m)$ for the total number of
FPL configurations with associated matching $X\cup m$, we have
\begin{equation} \label{eq:Summenformel} 
A_X(m)=
\sum _{\mathcal E,\mathcal F} ^{}c(\mathcal
E,\mathcal F)N(\mathcal E,\mathcal F,m,d),
\end{equation}
where the sum is over all possible choices of 
$\mathcal E\subseteq\{1,2,\dots,{d-2}\}$ and
$\mathcal F\subseteq\{1,2,\dots,{d-1}\}$ such that $\vert\mathcal
E\vert=\vert\mathcal F\vert-1$. In the next lemma, we record the
properties of the numbers
$N(\mathcal E,\mathcal F,m,d)$ which will allow us
to conclude the proof of Conjecture~\ref{conj:zub1} for $m< 3d$.

\begin{lem} \label{lem:cN}
For $m\ge 2d$, we have
\begin{enumerate} 
\item The number $N(\mathcal E,\mathcal F,m,d)$ is a polynomial in
$m$.
\item As a polynomial in $m$,
$m$ divides $N(\mathcal E,\mathcal F,m,d)$ for all $\mathcal E$
and $\mathcal F$, except if $\mathcal E=\{1,2,\dots,{d-2}\}$ and
$\mathcal F=\{1,2,\dots,{d-1}\}$, in which case
$N(\mathcal E,\mathcal F,m,d)=1$.
\end{enumerate}
\end{lem}

\begin{proof}
Our aim is to find a determinantal expression for $N(\mathcal
E,\mathcal F,m,d)$. To do so, we proceed as in the proof of
Theorem~\ref{main}, that is, we map the possible configurations in the
pentagonal region bijectively to rhombus tilings of a certain
region in the regular triangular lattice. 
As in the preceding proof, we draw triangles around free
vertices (where ``free'' has the same meaning as in that proof) in such
a way that two free vertices are neighbours
if and only if the corresponding triangles share an edge.
The region in the regular triangular lattice which we obtain for the
pentagonal region of Figure~\ref{fig:Ecke1} is shown in
Figure~\ref{fig:Region}. It is a hexagon with bottom side of length
$d-1$, lower-left side of length $d-1$, upper-left side of length $m-d+2$,
top side of length $d-2$, upper-right side of length $d$, and
lower-right side of length $m-d+1$. However, depending on the choice
of $\mathcal E$ and $\mathcal F$, along the top side and along the
upper-right side there are some notches (that is, triangles of unit
side length which are missing, as was the case for the
region $R(\mathcal E,m-3d+1)$ obtained in the proof of
Theorem~\ref{main}; compare with Figure~\ref{beau22}). 
In Figure~\ref{fig:Region},
the possible places for notches are labelled
$\{e_1,e_2,\dots,e_{d-2}\}$, respectively $\{f_1,f_2,\dots,f_{d-1}\}$.
The place labelled $f_d$ cannot be the place of a notch, and the
number of notches out of $\{f_1,f_2,\dots,f_{d-1}\}$ must be exactly
by 1 larger than the number of notches out of
$\{e_1,e_2,\dots,e_{d-2}\}$. An example of a choice of notches
for $d=5$, $m=7$, $\mathcal E=\{2\}$ and $\mathcal F=\{1,4\}$
(filled by a rhombus tiling of the resulting region) is shown
on the left in Figure~\ref{fig:region}.

\begin{figure}[h]
\psfrag{a1}{\large $e_1$}
\psfrag{a2}{\large $e_2$}
\psfrag{ad-1}{\large $e_{d-2}$}
\psfrag{b1}{\large $f_1$}
\psfrag{b2}{\large $f_2$}
\psfrag{bd+1}{\large $f_{d}$}
\psfrag{d}{\large $d-1$}
\psfrag{d+1}{\large $d$}
\psfrag{d-1}{\large $d-2$}
\psfrag{m-d}{\large $m-d+1$}
\psfrag{m-d+1}{\large $m-d+2$}
\resizebox*{10cm}{!}{\includegraphics{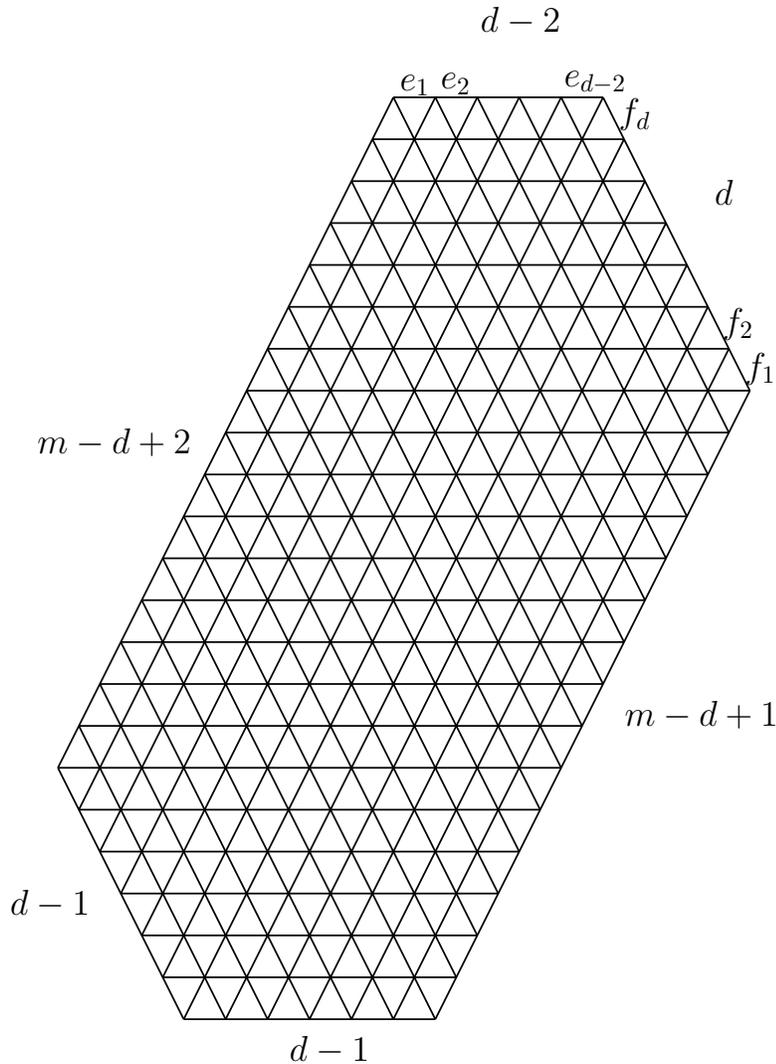}}
\caption{\label{fig:Region}The region to be tiled}
\end{figure}

Thus, the number $N(\mathcal
E,\mathcal F,m,d)$ is equal to the number of rhombus tilings of this
hexagonal region with notches. To get a formula for the number of
rhombus tilings, we apply the standard bijection between rhombus
tilings and families of non-intersecting lattice paths (see, e.g.,
\cite{CEKZ,CK}). 
The bijection is obtained as follows. One places vertices in each of
the mid-points of edges along the bottom-left side of the region,
and as well
in each mid-point along the downward edges which form part of a notch
on the top of the region, and along the downward edges along the
top-right side of the region. (In Figure~\ref{fig:region} these
midpoints are marked in boldface.) 
The vertices of the lower-left edges are subsequently connected to the
vertices on top and top-right by paths, by connecting the mid-points of
opposite downward edges in each rhombus of the tiling, see again
Figure~\ref{fig:region}.
Clearly, by construction, the paths are non-intersecting.
Subsequently, 
the paths are slightly rotated, and deformed so that they become
rectangular paths. Thus, one obtains families of paths with starting
points $A_i=(-i,i)$, $i=1,2,\dots,d-1$, and with 
end points a subset of cardinality $d-1$ from
the points $\{E_1,E_2,\dots,E_{d-2}\}\cup\{F_1,F_2,\dots,F_d\}$,
where $E_i=(-d+i,m+1)$, $i=1,\dots,d-2$, and $F_j=(d-j-1,m-d+j+1)$,
$j=1,2,\dots,d$. The point $F_d$ {\it must\/} be among the end points.
The family of paths which results
from our example rhombus tiling is shown on the right of
Figure~\ref{fig:region}. 

\begin{figure}[h]

$$
\hbox{\includegraphics{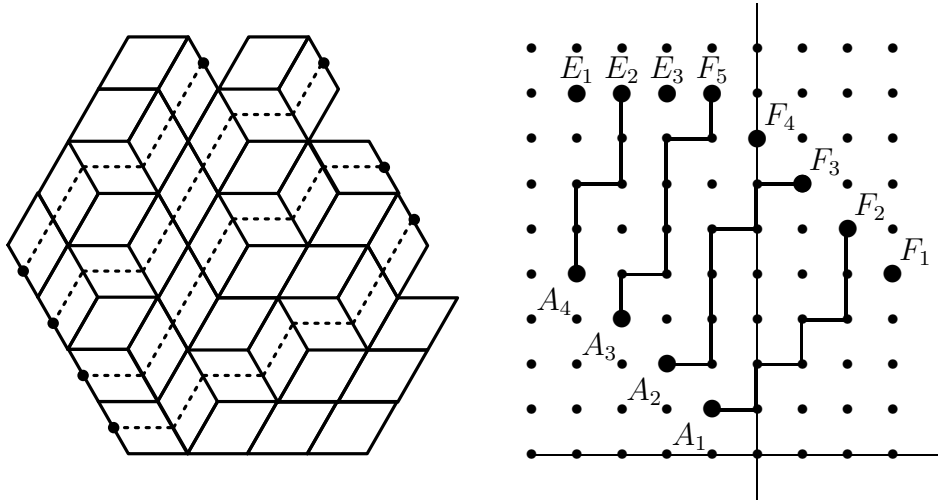}\kern4cm}
\Gitter(4,10)(-5,0)
\Koordinatenachsen(4,10)(-4,0)
\Pfad(-4,4),22122\endPfad
\Pfad(-3,3),2122212\endPfad
\Pfad(-2,2),1222121\endPfad
\Pfad(-1,1),1212122\endPfad
\DickPunkt(-1,1)
\DickPunkt(-2,2)
\DickPunkt(-3,3)
\DickPunkt(-4,4)
\DickPunkt(-4,8)
\DickPunkt(-3,8)
\DickPunkt(-2,8)
\DickPunkt(-1,8)
\DickPunkt(0,7)
\DickPunkt(1,6)
\DickPunkt(1,6)
\DickPunkt(2,5)
\DickPunkt(3,4)
\Label\lu{A_1}(-1,1)
\Label\lu{A_2}(-2,2)
\Label\lu{A_3}(-3,3)
\Label\lu{A_4}(-4,4)
\Label\ro{F_1}(3,4)
\Label\ro{F_2}(2,5)
\Label\ro{F_3}(1,6)
\Label\ro{F_4}(0,7)
\Label\o{\raise6pt\hbox{$F_5$}}(-1,8)
\Label\o{\raise6pt\hbox{$E_3$}}(-2,8)
\Label\o{\raise6pt\hbox{$E_2$}}(-3,8)
\Label\o{\raise6pt\hbox{$E_1$}}(-4,8)
\hskip2.5cm
$$
\caption{\label{fig:region}A rhombus tiling of the region with
notches,
and its corresponding family of non-intersecting lattice paths}
\end{figure}

Now we can apply Lemma~\ref{gv} to obtain a determinant for the number
of rhombus tilings. The number of paths from a starting point $A_i$ to an
end point, which forms an entry in the determinant in \eqref{eq:gv},
is equal to $\binom {m+1+j-d}{j+i-d}$ if the end point is $E_j$, and it is
$\binom {m}{d-j+i-1}$ if the end point is $F_j$. In both cases, this
is a polynomial in $m$, and thus claim~(1) follows. 

To prove claim~(2), we observe that $m$ divides
the second of the above binomial coefficients, $\binom {m}{d-j+i-1}$,
as long as $d-j+i-1>0$. However, this will be always the case, except
if $i=1$ and $j=d$. In particular, if we choose an $F_j$ with $j<d$ as
an end point, all the entries in the column corresponding to $F_j$ 
in the determinant \eqref{eq:gv} will be divisible by $m$, and, hence,
the determinant itself. Thus, the only case where the determinant is
not divisible by $m$ is the one where the set of end points is
$\{E_1,E_2,\dots,E_{d-2},F_d\}$. This corresponds to the choice
$\mathcal E=\{1,2,\dots,{d-2}\}$ and
$\mathcal F=\{1,2,\dots,{d-1}\}$. In addition, in that case 
the determinant in \eqref{eq:gv} is equal to 1 because it is the
determinant of a triangular matrix with 1s on the antidiagonal. 
This completes the proof of the lemma.
\end{proof}

We are now in the position to complete the proof of
Conjecture~\ref{conj:zub1} for $m<3d$. We wish to prove that the
polynomial which results from the right-hand side of 
\eqref{eq:Summenformel} by substituting
the determinantal formula for $N(\mathcal E,\mathcal F,m,d)$ 
obtained in the preceding proof of
Lemma~\ref{lem:cN} gives the number of FPL configurations under
consideration {\it for all $m\ge0$}.
As we remarked earlier,
the arguments so far show only that this is indeed the case for $m\ge2d$.

We verify next that this is also the case for $m=0$. If we put $m=0$
in the right-hand side of \eqref{eq:Summenformel} 
(with the afore-mentioned substitution for $N(\mathcal E,\mathcal F, m,d)$),
then, by
Lemma~\ref{lem:cN}, the only term which survives is the one for
$\mathcal E=\{1,2,\dots,{d-2}\}$ and
$\mathcal F=\{1,2,\dots,{d-1}\}$. In addition, in that case
we have $N(\mathcal E,\mathcal F, m,d)=1$. Thus, if $m=0$, the expression in
\eqref{eq:Summenformel} is equal to $c(\mathcal E,\mathcal F)$,
the number of possible configurations in the triangular region in the
upper-right corner of the square grid. Figure~\ref{fig:Ecke2} shows
what happens inside this triangular region for this choice of
$\mathcal E$ and $\mathcal F$: from the external links occupied by the
matching $X$ there propagate zig-zag lines of fixed edges into the
interior, so that only a square region with side length $d-1$ remains
undetermined. Thus, the
number of possible FPL configurations inside this triangular region is
indeed equal to the number of FPL configurations with
associated matching $X$ (that is, FPL configurations {\it on the
square $Q_d$} with associated matching $X$), 
which is exactly what we wanted to prove.

While it seems that we have still a large gap (namely the values of
$m$ between $1$ and $2d-1$) to overcome, the assertion now follows:
let $H_X(m)$ denote the polynomial on the right-hand side of
\eqref{eq:Summenformel}.
By the above arguments we know that,
for an arbitrary non-negative integer $m$, 
the number $A_X(m)=A_{X\cup m}(0)$ 
of FPL configurations with associated matching
$X\cup m$ is equal to $H_{X\cup m}(0)$.
Furthermore, for any non-negative integer $s$ we have
$A_{X\cup m}(s)=A_X(m+s)$. By the preceding
arguments, if $s$ is sufficiently large, we have also $A_{X\cup m}(s)=
H_{X\cup m}(s)$ and $A_{X}(m+s)=H_{X}(m+s)$.
However, it also follows from the preceding arguments that 
$H_{X\cup m}(s)$ and $H_{X}(m+s)$ are polynomials in $s$.
Since they agree for an infinite number of values of $s$, they must be
identical. In particular, $A_{X\cup m}(0)=H_{X\cup m}(0)=
H_{X}(m)$. Thus, the number $A_{X\cup m}(0)$ 
of FPL configurations with associated
matching $X\cup m$ is indeed given by the same polynomial for any
$m$. It must necessarily be equal to the polynomial found in
Section~\ref{secCon} (see the last paragraph of that section). 
Conjecture~\ref{conj:zub1} is now completely proved.

\section{\label{secGen}Proof of Conjecture~\ref{conj:zub2} for $m$
large enough}

In this section we show how the ideas developed in the proof of 
Theorem~\ref{main} can be extended to prove Conjecture~\ref{conj:zub2}
for large enough $m$.

\begin{figure}
\psfrag{T}{$X$}\psfrag{S}{$Y$}
\includegraphics{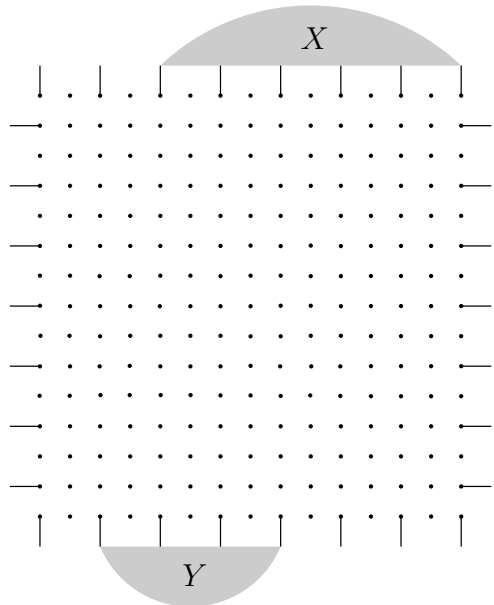}

\caption{\label{beau23}Placing the matching around the grid}
\end{figure}

Let $X$ and $Y$ be two non-crossing matchings with $d$ and $e$ arches,
respectively. Without loss of generality, we may suppose that $d\geq
e$. We choose to place the matching
$X\cup m\cup Y$ in Conjecture~\ref{conj:zub2} 
in such a way that, again, the set of arches
of $X$ appear on the very right of the upper side of
the grid $Q_{n}=Q_{d+e+m}$, see Figure~\ref{beau23} for a schematic picture,
that is, we place these arches on the external links labelled 
$ n-4d+2, n-4d+4,\ldots, n-2, n$. In order to guarantee that $X$ has
place along the upper side of the square grid, we must assume that
$m\ge 3d-e$.

\begin{figure}
\psfrag{T}{$X$}\psfrag{S}{$Y$}
 \psfrag{A}{$A$}
 \psfrag{B}{$B$}
 \psfrag{C}{$C$}
 \psfrag{D}{$D$}
 \psfrag{E}{$E$}
 \psfrag{J}{$J$}
 \psfrag{K}{$K$}
 \psfrag{L}{$L$}
 \psfrag{F}{$F$}
 \psfrag{G}{$G$}
 \psfrag{M}{$M$}
\psfrag{rho}{$\xi_1$}\psfrag{sig}{$\eta_1$}
\psfrag{rho'}{$\xi_2$}\psfrag{sig'}{$\eta_2$}
\resizebox*{10cm}{!}{\includegraphics{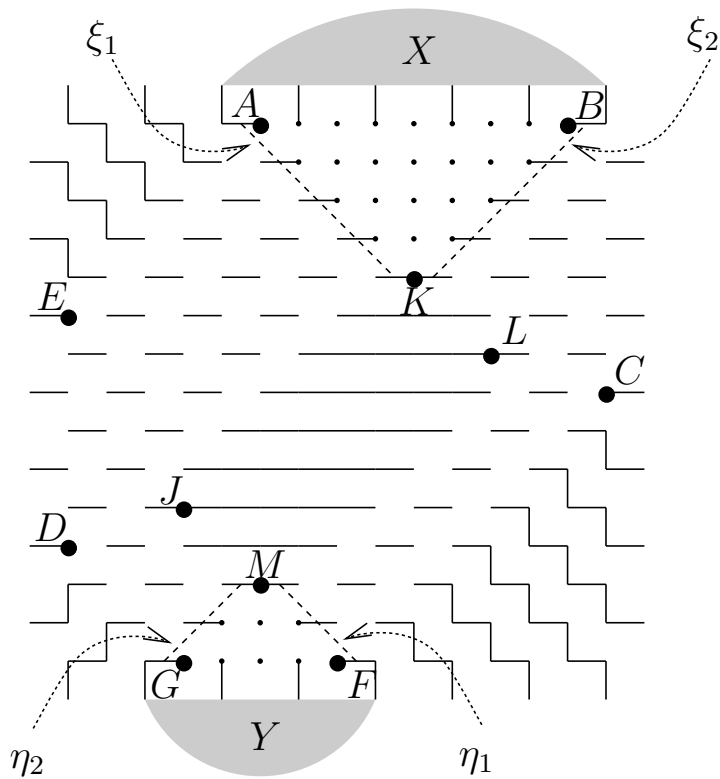}}

\caption{\label{beau24}The fixed edges}
\end{figure}

Next we determine the set of fixed edges using Lemma~\ref{fixedg}.
For convenience (the reader should consult Figure~\ref{beau24} while
reading the following definitions; it contains an example with $n=15$,
$d=3$, $e=2$ and $m=10$), we let $A,B,C,F,G,D,E$ (in this
order!) be the border vertices 
of the external links labelled $n-4d+3$, $n-1$, $n+2d+2e-2$, $-n-2d-2e+3$,
$-n-2d+2e-1$, $-n+2d-2e+2$, $-n+4d-2$, respectively, 
we let $K$ be the point in the interior of $Q_n$ which makes $ABK$
into a rectangular isosceles triangle, with the right angle at $K$,
we let
$M$ be the analogous point which makes $FGM$
into a rectangular isosceles triangle, with the right angle at $M$, 
we let $J$ be
the intersection point of the line connecting $F$ and $M$ and
the line connecting $B$ and $K$, and we let $L$ be
the intersection point of the line connecting $G$ and $M$ and
the line connecting $A$ and $K$.
We state the result of the application of Lemma~\ref{fixedg} to the
current case in form of the following lemma.

\begin{lem} \label{lem:fixed2}
The region of fixed edges of the FPL configurations corresponding
to the matching $X\cup m\cup Y$ contains all the
edges indicated in Figure~\emph{\ref{beau24}}, that is:

\begin{enumerate} 
\item all the horizontal edges in the rectangular region $JKLM$,
\item every other horizontal edge in the L-shaped region $AKJMGDE$ as
indicated in the figure,
\item every other horizontal edge in the region $BCFMLK$ as
indicated in the figure,
\item the zig-zag lines in the corner regions above the line $EA$,
respectively below the lines $DG$ and $FC$, as indicated in the figure.
\end{enumerate}
\end{lem}

Let $\xi_1$ be the segment which
connects the point which is half a unit to the left of $A$ and the point
which is half a unit to the left of $K$, and let
$\xi_2$ be the segment which
connects the point which is half a unit to the right of $K$ 
with the point which is half a unit to 
the right of $B$, see Figure~\ref{beau24}.
Similarly, let $\eta_1$ be the segment which
connects the point which is half a unit to the left of $G$ and the point
which is half a unit to the left of $M$, and let
$\eta_2$ be the segment which
connects the point which is half a unit to the right of $M$ 
with the point which is half a unit to 
the right of $F$, see again Figure~\ref{beau24}.
There are $2d-2$ vertical edges which cross $\xi_1$, and the same is true for
$\xi_2$, and there are $2e-2$ vertical edges which cross $\et_1$, and
the same is true for $\et_2$.

The next lemma says that, also in this more general situation, it is
exactly one half of the vertical edges which cross $\xi_1$, $\xi_2$, $\et_1$,
and $\et_2$, respectively, 
which are taken by any FPL configuration with associated
matching $X\cup m\cup Y$.

\begin{lem} \label{lem:de}
Any FPL configuration with associated matching $X\cup m\cup Y$,
occupies exactly $d-1$ vertical edges crossing the segments $\xi_{1}$,
the same being true for $\xi_2$,
and it occupies exactly $e-1$ vertical edges crossing the segments 
$\eta_{1}$, the same being true for $\eta_2$.
\end{lem}

\begin{proof}
This follows from Corollary~\ref{cor}.(2).
\end{proof}

Hence, any such choice of edges as described in Lemma~\ref{lem:de}
can be encoded, as before, 
by four sets $\mathcal E_1, \mathcal E_2,\mathcal
F_1,\mathcal F_2$, the sets $\mathcal E_1$ and $\mathcal E_2$ being subsets of 
$\{1,2,\dots,2d-2\}$, and the sets
$\mathcal F_1$ and $\mathcal F_2$ being subsets of
$\{1,2,\dots,2e-2\}$. We make the convention
that $\mathcal E_1$ encodes the set of vertical edges crossing 
the segment $\xi_{1}$, $\mathcal E_2$ encodes the set of vertical 
edges crossing the segment $\xi_{2}$,
$\mathcal F_1$ the set of edges crossing the segment $\eta_{1}$ and 
finally $\mathcal F_2$
the set of edges crossing the segment $\eta_{2}$. 
We denote by $a_{X}(\mathcal E_1,\mathcal E_2)$
the number of FPL configurations in the triangular region $ABK$, 
which ``obey'' the matching $X$, which, out of the edges crossing the 
segment $\xi_{1}$, occupy exactly those
determined by $\mathcal E_1$, and which, out of the edges crossing 
the segment $\xi_{2}$, occupy exactly those
determined by $\mathcal E_2$. Note that, by Corollary~\ref{cor}.(1),
$a_{X}(\mathcal E_1,\mathcal E_2)$ can only be non-zero
if $\la(\mathcal E_1)\subseteq \la(X)$ and $\la(\mathcal E_2)\subseteq
\la(X)'$, 
where $\la(X)'$ denotes the Ferrers diagram conjugate to $\la(X)$. 
For the latter assertion, we use the following trivial observation:
Let $X$ be a non-crossing matching with $d$ arches, 
and let $X^r$ be the ``reversed''
matching, that is, whenever $i$ is matched with $j$ in $X$, $2d-i+1$ is
matched with $2d-j+1$ in $X^r$. Then, under
the correspondence between matchings and Ferrers diagrams described in
Section~\ref{sec:Ferrers}, we have $\la(X^r)=\la(X)'$.

\begin{figure}
\psfrag{h}{$h$}\psfrag{k}{$k$}
\psfrag{d+e}{$\kern-7ptd+e$}\psfrag{0}{$0$}
\psfrag{1}{$1$}\psfrag{d-code of X}{$d$-code of $X$}
\psfrag{e-code of Y}{$e$-code of $Y$}
\includegraphics{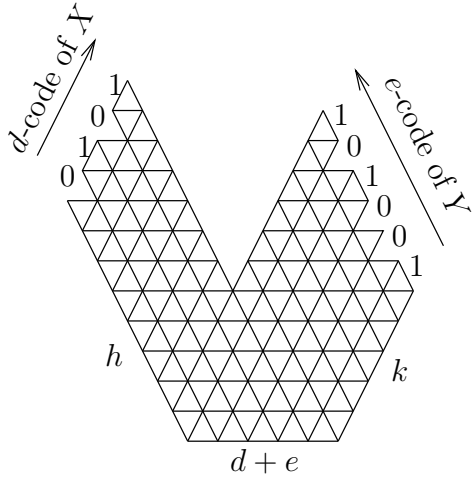}

\caption{\label{beau25}The region $R(X,Y,d,e,h,k)$}
\end{figure}

In order to proceed, 
we need to introduce a region in the triangular lattice
which is parametrized by two Ferrers diagram and four non-negative integers.
Let $h,k,d,e$ be non-negative integers, and let 
$\la\subseteq(d^{d})$
and $\mu\subseteq(e^{e})$ 
be two Ferrers diagrams. Then we define
$R(X,Y,d,e,h,k)$ to be the V-shaped region (the reader should consult
Figure~\ref{beau25} while reading the following description; the
figure illustrates $R(X,Y,d,e,h,k)$ 
for $d=2$, $e=3$, $h=8$, $k=5$, $\la=(1,0)$, and
$\mu=(2,1,1)$) with base side of length $d+e$, left side $h$,
followed by a side in direction north-east of length $2d$ with notches 
which will be explained in just a moment, a V-shaped ``valley,''
a side in direction south-east of length $2e$ with notches 
which will be explained in just a moment, and finally a right side of
length $k$. To determine the notches along the side in direction
north-east, we read the $d$-code of $\la$ and while, at the same time,
moving along the side from south-west to north-east, we
put a notch whenever we read a $0$, and we leave a horizontal
piece whenever we read a $1$. The notches along the side in direction
south-east are determined in a similar way from the $e$-code of $\mu$,
however, while reading the $e$-code, we move
along the side from south-east to north-west. 
Let $r(X,Y,d,e,h,k)$ be the number of rhombus tilings
of $R(X,Y,d,e,h,k)$. 

With the above notation, we have the following expansion for the
number $A_{X,Y}(m)$ of FPL configurations with associated matching
$X\cup m\cup Y$.

\begin{thm} \label{main2}
Let $X$ and $Y$ be two non-crossing matchings with $d$ and $e$ arches,
respectively. Then, for $m\geq3d-e$, 
\begin{multline} \label{eq:main2}
A_{X,Y}(m)=\sum_{\mathcal E_1,\mathcal E_2,\mathcal F_1,\mathcal F_2}
a_{X}(\mathcal E_1,\mathcal E_2)\,a_{Y}(\mathcal F_1,\mathcal F_2) \\
   \cdot r(\la(\mathcal E_1),\la(\mathcal F_2),d-1,e-1,m-3d+e+2,2d-2e+1)\\
\cdot
 r(\la(\mathcal F_1),\la(\mathcal E_2),e-1,d-1,m-d-e+2,1),
\end{multline}
where the sum is taken over all subsets $\mathcal E_1,\mathcal E_2$
of $\{1,2,\dots,2d-2\}$ and all subsets $\mathcal F_1,\mathcal F_2$
of $\{1,2,\dots,2e-2\}$ such that 
$\la(\mathcal E_1)\subseteq \la(X)$,
$\la(\mathcal E_2)\subseteq \la(X)'$, 
$\la(\mathcal F_1)\subseteq \la(Y)$, and 
$\la(\mathcal F_2)\subseteq \la(Y)'$.
\end{thm}

\begin{proof}
Let us fix $d-1$ edges which cross $\xi_1$, 
$d-1$ edges which cross $\xi_2$, 
$e-1$ edges which cross $\et_1$, 
$e-1$ edges which cross $\et_2$, encoded by the sets
$\mathcal E_1,\mathcal E_2\subseteq\{1,2,\dots,2d-2\}$ of cardinality
$d-1$ and by the sets
$\mathcal F_1,\mathcal F_2\subseteq\{1,2,\dots,2e-2\}$ of cardinality
$e-1$, respectively. Clearly, the number of FPL configurations with
associated matching $X\cup m\cup Y$ which occupy exactly these
vertical edges crossing $\xi_1$, $\xi_2$, $\et_1$, and $\et_2$
decomposes into a product: it equals the number of configurations in
the region $ABK$, given by $a_{X}(\mathcal E_1,\mathcal E_2)$, times
the number of configurations in
the region $FGM$, given by $a_{Y}(\mathcal F_1,\mathcal F_2)$, times
the number of configurations in
the region $AKJMGDE$, times
the number of configurations in
the region $FMLKBC$. 

In order to compute the number of configurations in the latter two
regions, as in the proof of Theorem~\ref{main}, 
we translate again the problem of enumerating FPL
configurations into a problem of enumerating rhombus tilings. 
If we do this for the region $AKJMGDE$, then, as a result, we have to
count all the rhombus tilings of the region $R(\la(\mathcal
E_1),\la(\mathcal F_2),d-1,e-1,m-3d+e+2,2d-2e+1)$, while for the
region $FMLKBC$ we arrive at the problem of counting all the rhombus
tilings of the region $R(\la(\mathcal F_1),\la(\mathcal
E_2),e-1,d-1,m-d-e+2,1)$. This finishes the proof of the theorem.
\end{proof}

Our aim is to show that the quantities $r(.)$
appearing in \eqref{eq:main2} are
polynomials in $m$, and that the ``dominating'' term
in the sum in \eqref{eq:main2} comes from
choosing $\mathcal E_1,\mathcal E_2,\mathcal F_1,\mathcal F_2,$
such that $\la(\mathcal E_1)=\la(X)$, $\la(\mathcal F_1)=\la(Y)$, and 
$\la(\mathcal E_2)=\la(\mathcal F_2)=\emptyset$.
In order to do so, we need several preparatory lemmas.

\begin{lem} \label{lem:1}
Let $\la$ be a partition contained in $(d^d)$, and let $\mu$ be a
partition contained in $(e^e)$. Then
the number $r(\la,\mu,d,e,h,k)$ of rhombus tilings of the region
$R(\la,\mu,d,e,h,k)$ is equal to
\[
r(\la,\mu,d,e,h,k)=\frac{1}{|\la|!}p_{\la,\mu,d,e,k}(h),\]
where $p_{\la,\mu,d,e,k}(h)$ is a polynomial of degree $|\la|$ with integer
coefficients. 
\end{lem}
\begin{proof}
Again, we use the standard bijection between rhombus tilings and
non-intersecting lattice paths (cf.\ the proof of
Lemma~\ref{lem:cN}). This time, we introduce vertices along the bottom
side of the region $R(\la,\mu,d,e,h,k)$, and along the horizontal
edges which are part of the notches on top-left and top-right of the
region (see Figure~\ref{beau33}, which shows the region
$R((1,0),(2,1,1),2,3,8,5)$ of Figure~\ref{beau25}). These vertices are
then connected by paths, as before. If, finally, the paths are again
rotated and slightly deformed so that they become rectangular paths,
we obtain families $(P_1,P_2,\dots,P_{d+e})$ of non-intersecting
lattice paths, where the starting points of the paths $P_i$ are
$A_i=(-i,i)$, $i=1,2,\dots,d+e$, and where their end points are the points
$E_j=(\la_j-j-e,h+d+e)$, $j=1,2,\dots,d$, and $F_j=(k-1,\mu_j+e-j+1)$,
$j=1,2,\dots,e$.

\begin{figure}
$$
\hbox{\includegraphics{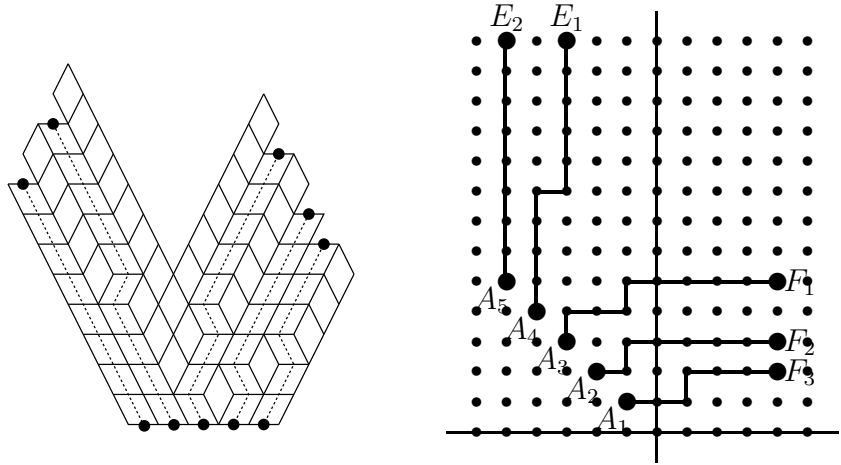}\kern4cm}
\Einheit.4cm
\Gitter(6,14)(-6,0)
\Koordinatenachsen(6,14)(-6,0)
\Pfad(-1,1),112111\endPfad
\Pfad(-2,2),1211111\endPfad
\Pfad(-3,3),211211111\endPfad
\Pfad(-4,4),2222122222\endPfad
\Pfad(-5,5),22222222\endPfad
\DickPunkt(-1,1)
\DickPunkt(-2,2)
\DickPunkt(-3,3)
\DickPunkt(-4,4)
\DickPunkt(-5,5)
\DickPunkt(4,2)
\DickPunkt(4,3)
\DickPunkt(4,5)
\DickPunkt(-3,13)
\DickPunkt(-5,13)
\Label\lu{A_1}(-1,1)
\Label\lu{A_2}(-2,2)
\Label\lu{A_3}(-3,3)
\Label\lu{A_4}(-4,4)
\Label\lu{A_5}(-5,5)
\Label\o{\raise6pt\hbox{$E_2$}}(-5,13)
\Label\o{\raise6pt\hbox{$E_1$}}(-3,13)
\Label\r{F_3}(4,2)
\Label\r{F_2}(4,3)
\Label\r{F_1}(4,5)
\hskip2cm
$$
\caption{Non-intersecting lattice paths for 
a rhombus tiling of $R(\la,\mu,d,e,h,k)$\label{beau33}}
\end{figure}

We may now apply Lemma~\ref{gv} to obtain a determinant for the
number these families of non-intersecting lattice paths. The number of
paths from a starting point $A_i$ to an end point, which forms an entry in the
determinant in \eqref{eq:gv}, is equal to $\binom
{\la_j-j+h+d}{\la_j-j-e+i}$ if the end point is $E_j$,
and it is $\binom {k+\mu_j+e-j}{k-1+i}$ if the end point is $F_j$. 

If we expand the determinant \eqref{eq:gv}, with these
specializations, according to its
definition, then we obtain, up to an overall sign,
\begin{equation} \label{eq:detgv}
\sum _{\si\in S_{d+e}} ^{}\sgn\si 
\Bigg(
\prod _{j=1} ^{d}\binom {\la_j-j+h+d}{\la_j-j-e+\si(j)}\Bigg)
\Bigg(\prod _{j=d+1} ^{d+e} \binom {k+\mu_j+e-j}{k-1+\si(j)}\Bigg),
\end{equation}
where $S_N$ denotes the group of permutations of $\{1,2,\dots,N\}$.

To determine the degree of \eqref{eq:detgv} as a polynomial in $h$,
we observe that the second product does not contribute, whereas the
degree in $h$ of the first product is $\sum _{j=1}
^{d}(\la_j-j-e+\si(j))=\vert\la\vert+\sum _{j=1}
^{d}(\si(j)-j-e)$. The latter sum is maximal exactly if $\si(j)=j+e$,
$j=1,2,\dots,d$, in which case it equal to 0. Thus, the degree of
\eqref{eq:detgv} as a polynomial in $h$ is indeed $\vert\la\vert$. 

Finally, to show that, if we put \eqref{eq:detgv} on the common
denominator $\vert\la\vert!$, the numerator polynomial has integer
coefficients, we observe that the second product in \eqref{eq:detgv}
is an integer. Thus, it suffices to multiply the first product by
$\vert\la\vert!$, and show that the result is a polynomial in $h$ with
integer coefficients. Indeed, this is
\begin{multline} \label{eq:Nenner}
\vert\la\vert!\Bigg(
\prod _{j=1} ^{d}\binom {\la_j-j+h+d}{\la_j-j-e+\si(j)}\Bigg)\\=
L!\binom {\vert\la\vert}{\la_1-1-e+\si(1),\dots,\la_d-d-e+\si(d),L}
\prod _{j=1} ^{d}(h+d+e-\si(j)+1)_{\la_j-j-e+\si(j)},
\end{multline}
where $L=
\sum _{j=1} ^{d}(j+e-\si(j))$ (which is non-negative, as we have just
shown!), where the Pochhammer symbol $(\al)_s$ is defined by
$(\al)_s:=\al(\al+1)\cdots(\al+s-1)$, and where for
$N=N_1+N_2+\dots+N_s$ the multinomial
coefficient is defined by
$$\binom N{N_1,N_2\dots,N_s}=\frac {N!} {N_1!\,N_2!\cdots N_s!}.$$
Evidently, the right-hand side in \eqref{eq:Nenner} is a polynomial in $h$
with integer coefficients.
\end{proof}

\begin{lem} \label{lem:2}
Let $\la$ be a partition contained in $(d^{d})$.
The number of rhombus tilings of $R(\la,\emptyset,d,e,h,k)$ is equal
to the number of rhombus tilings of the region $R(\la,d,h)$
defined in Section~\ref{secHook}.
\end{lem}
\begin{figure}[h]
\begin{center}\includegraphics{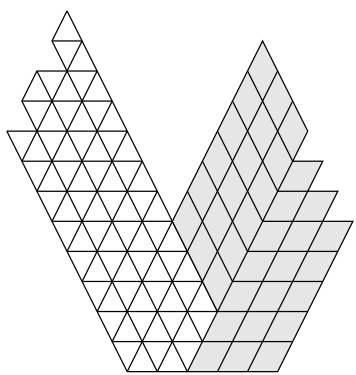}\end{center}

\caption{The forced rhombi in $R(\la,\emptyset,d,e,h,k)$ \label{beau34}}
\end{figure}
\begin{proof}
The reader should consult Figure~\ref{beau34} for an example of a
region $R(\la,\mu,d,e,h,k)$ where $\mu$ is the empty partition.
In that case, there are many ``forced'' rhombi in the right half of
the V-shaped region, that is, rhombi which are common to every rhombus
tiling of the region. In the example in Figure~\ref{beau34}, these
rhombi are the ones in the shaded region. Once we remove the forced
rhombi, we are left with (a reflected version of) 
the region $R(\la,d,h)$.
\end{proof}

\begin{lem} \label{lem:3}
Let $X$ be a non-crossing matching with $d$ arches, 
and let $\mathcal E_1$ be the subset
of $\{1,2,\dots,2d-2\}$ for which $\la(X)=\la(\mathcal E_1)$. Then \[
a_{X}(\mathcal E_1,\mathcal E_2)=\left\{ \begin{array}{ll}
0, & \textrm{if $\mathcal E_2\neq\{1,2,\dots,d-1\}$,}\\
1, & \textrm{if $\mathcal E_2=\{1,2,\dots,d-1\}.$}\end{array}\right.\]
\end{lem}

\begin{proof}This follows from the proof of Corollary~\ref{cor}.(3).
\end{proof}

\begin{remark}
The reader should note that the choice of $\mathcal
E_2=\{1,2,\dots,d-1\}$ means that $\la(\mathcal E_2)=\emptyset$.
\end{remark}

If we use Lemmas~\ref{lem:1}--\ref{lem:3} in Theorem~\ref{main2}, then
we obtain the second main result of this paper.

\begin{thm} \label{main2a}
Let $X$ and $Y$ be two non-crossing matchings with $d$ and $e$ arches,
respectively. Then, for $m\geq3d-e$,
\begin{multline} \label{eq:main2a}
A_{X,Y}(m)  =  SSYT(\la(X),m-2d+e+1)\cdot SSYT(\la(Y),m-d+1) \\
+\sum_{\mathcal E_1,\mathcal E_2,\mathcal F_1,\mathcal
F_2}\Big(a_{X}(\mathcal E_1,\mathcal E_2)\,a_{Y}(\mathcal F_1,\mathcal
F_2)\kern5.5cm\\
 \cdot r(\la(\mathcal E_1),\la(\mathcal F_2),d-1,e-1,m-3d+e+2,2d-2e+1)\\
\cdot r(\la(\mathcal
F_1),\la(\mathcal E_2),e-1,d-1,m-d-e+2,1)\Big), 
\end{multline}
where the sum is taken over all subsets $\mathcal E_1,\mathcal E_2$
of $\{1,2,\dots,2d-2\}$ and all subsets $\mathcal F_1,\mathcal F_2$
of $\{1,2,\dots,2e-2\}$ such that 
$\la(\mathcal E_1)\varsubsetneq \la(X)$,
$\la(\mathcal E_2)\subseteq \la(X)'$, 
$\la(\mathcal F_1)\varsubsetneq \la(Y)$, and 
$\la(\mathcal F_2)\subseteq \la(Y)'$.
\end{thm}

\begin{proof}
In the sum on the right-hand side of \eqref{eq:main2} we single out
the term where $\la(\mathcal E_1)=\la(X)$, $\la(\mathcal
E_2)=\la(\mathcal F_2)=\emptyset$, and $\la(\mathcal F_1)=\la(Y)$. By
Lemma~\ref{lem:3}, in that case we have
$a_X(\mathcal E_1,\mathcal E_2)=a_Y(\mathcal F_1,\mathcal F_2)=1$, and by
Lemma~\ref{lem:2} we have
$$r(\la(X),\emptyset,d-1,e-1,m-3d+e+2,2d-2e+1)=r(\la(X),d-1,m-3d+e+2)$$
and 
$$ r(\la(Y),\emptyset,e-1,d-1,m-d-e+2,1)=r(\la(Y),e-1,m-d-e+2),$$
where $r(\la,d,h)$ is short for the number of rhombus tilings of the
region $R(\la,d,h)$. By Theorem~\ref{thm:rhla}, the latter numbers are
given by the corresponding specializations in formula
\eqref{eq:hook-content}. Thus, we obtain the first term on the
right-hand side of \eqref{eq:main2a}. The asserted description of the
summation range of the sum on the right-hand side of \eqref{eq:main2a}
follows then from Lemma~\ref{lem:3}.
\end{proof}

Zuber's Conjecture~\ref{conj:zub2}, in the case that $m\ge 3d-e$, 
is now a simple corollary of the above theorem.

\begin{proof}[Proof of Conjecture~\ref{conj:zub2} for $m\ge 3d-e$]
The polynomiality in $m$ of $A_{X,Y}(m)$ is obvious from \eqref{eq:main2} and
Lemma~\ref{lem:1}. The assertion about the integrality of the
coefficients of the ``numerator'' polynomial $P_{X,Y}(m)$ follows as well
from Lemma~\ref{lem:1}. Finally,
to see that the leading coefficient of $P_{X,Y}(m)$ is 
$\dim(\la(X))\cdot\dim(\la(Y))$, 
one first observes that the leading term in
\eqref{eq:main2a} appears in the first term on the right-hand side.
The claim follows now by a
combination of \eqref{eq:hook-content} and \eqref{eq:hook}.
\end{proof}

\begin{remark}
The reader may wonder why we are not able to prove
Conjecture~\ref{conj:zub2} for the ``small'' values of $m$, that is,
for the range $0\le m<3d-e$. While we could try to imitate the
approach of Section~\ref{secCon1}, and place $X$ around the top-right
corner and $Y$ around the bottom-left corner of the square grid $Q_n$,
the problem arises when
Lemma~\ref{fixedg} is applied to determine the fixed
edges. Unfortunately, here, it is not the case that every grid vertex
is on at least one fixed edge. On the contrary, there will be a strip
without any fixed edge extending between the places where $X$ and $Y$
sit. Hence, the rhombus tiling arguments of Section~\ref{secCon1} do
not apply.
\end{remark}

\section{Auxiliary results for FPL configurations in a triangle}
\label{sec:aux}

The purpose of this section is to provide the proofs of
Lemmas~\ref{imposs}, \ref{lem:de} and \ref{lem:3}. These will be
consequences of Theorem~\ref{th:main} below.

\begin{figure}[h]
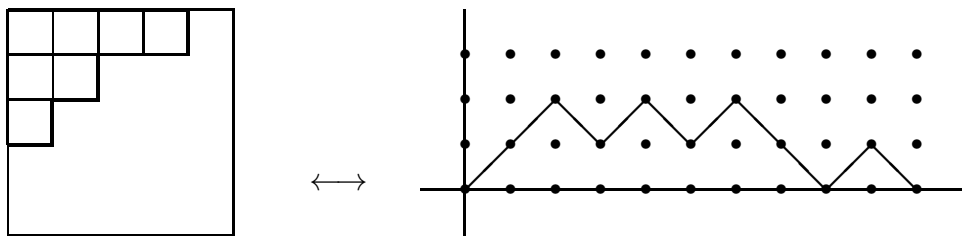

$$ 
\Pfad(0,1),12121125555666\endPfad
\PfadDicke{.5pt}
\Pfad(0,-1),11111222225555566666\endPfad
\Pfad(0,2),122\endPfad
\Pfad(0,3),112\endPfad
\Pfad(3,3),2\endPfad
\hbox{\hskip4cm$\longleftrightarrow$\hskip1.3cm}
\Gitter(11,4)(0,0)
\Koordinatenachsen(11,4)(0,0)
\PfadDicke{1pt}
\Pfad(0,0),3343434434\endPfad
\hskip5.5cm
$$
\caption{A Ferrers diagram embedded in a square and its corresponding
path \label{fig:Pfad}}
\end{figure}

For convenience, in this section we will encode Ferrers diagrams $\la$
which come from non-crossing matchings $X$ or sets $\mathcal E$
in form of paths consisting of up-steps $(1,1)$
and down-steps $(1,-1)$, which start at the origin. To obtain the path
corresponding to a Ferrers diagram $\la$ contained in the square
$(d^d)$, we trace the lower-right contour of $\la$
(including the parts along the left side and the upper side of the
square), and rotate it by
$45^\circ$. See Figure~\ref{fig:Pfad} for an illustration in the case
of $\la=(4,2,1)$ when viewed as a Ferrers diagram contained in $(5^5)$.

Given a non-crossing matching $X$ with $d$ arches,
we shall write $\pi(X)$ for the path corresponding to $\la(X)$ as
contained in $(d^d)$. For example, for $X_0$ the non-crossing matching
$\{1\leftrightarrow 8,2\leftrightarrow 3,4\leftrightarrow 5,
6\leftrightarrow 7,9\leftrightarrow 10\}$, the path $\pi(X_0)$ is the
one in Figure~\ref{fig:Pfad}.
Clearly, $\pi(X)$ will always be a Dyck path, that is, a path which
never passes below the $x$-axis and finally ends on the $x$-axis. 

\begin{figure}[htb]
\psfrag{A}{\small$A$}
\psfrag{B}{\small$B$}
\psfrag{K}{\small$K$}
\psfrag{X}{\small$X$}               
\psfrag{x1}{\small$\xi_1$}
\psfrag{PE}{\small$\pi({\mathcal E})$}
\psfrag{PX}{\small$\pi(X)$}
\centering \includegraphics[width=1\textwidth]{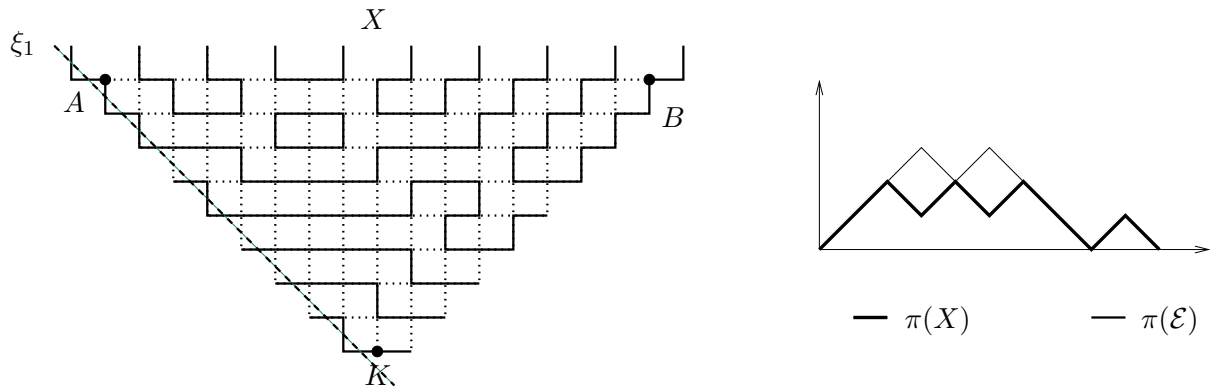}
\caption{The triangular region $\mathcal{R}$ with an example
of an FPL configuration, and the paths $\pi({\mathcal E})$ 
and $\pi(X)$.} \label{exemple}
\end{figure}

On the other hand, given a subset $\mathcal
E$ of $\{1,2,\ldots,{2d-2}\}$, let $c_{\mathcal E}$ be the
corresponding binary string as described at the beginning of
Section~\ref{secCon}. 
We shall denote by $\pi(\mathcal E)$
the path which arises from the string $0c_{\mathcal E}1$ (i.e., we
prepended a 0 and appended a 1 to $c_{\mathcal E}$) by starting the
path at the origin, and then reading $0c_{\mathcal E}1$ from left to
right, interpreting a 0 as an up-step $(1,1)$, and a 1 as a down-step
$(1,-1)$. For example, if $d=5$ and $\mathcal E_0=\{1,2,4,8\}$, then
$c_{\mathcal E_0}=00101110$, and the path $\pi(\mathcal E_0)$ is the thin path
on the right in Figure~\ref{exemple}.
It should be noted that, in general, $\pi(\mathcal E)$ need
not at all be a Dyck path; neither does it have to end on the
$x$-axis, nor does it have to stay above the $x$-axis.

Now let us consider a triangular region $\mathcal R$ 
below a non-crossing matching
$X$ with $d$ arches, 
as the region $ABK$ in Figures~\ref{beau15} or \ref{beau24}.
(The same arguments will apply to the region $FGM$ in
Figure~\ref{beau24}.) We also include the $2\times(2d-1)$ vertices on the
left and the right.
See the left part of Figure~\ref{exemple}
for an example with $d=5$. 
As earlier, we denote the segment connecting the points half a unit to
the left of $A$ and $K$ by $\xi_1$. In accordance with the situations
which we face in Sections~\ref{secCon} and \ref{secGen}, the FPL
configurations that we consider in this section will always contain all
the horizontal edges along the left border of $\mathcal R$,
that is, the horizontal edges crossing the segment $\xi_1$, and all
the horizontal edges along the right border of $\mathcal R$.
Moreover, we will only consider FPL configurations consisting
exclusively of loops connecting external edges occupied by $X$ or
loops which enter $\mathcal R$ from the left and exit $\mathcal R$ on
the right. (The latter correspond to some of the $m$ ``parallel''
loops in our two problems.) Let us call these FPL configurations 
{\it $\mathcal R$-FPL configurations}.

Given an $\mathcal R$-FPL configurations, we encode the set
of vertical edges crossing $\xi_1$ which are taken by the configuration 
by a subset $\mathcal E$ of
$\{1,2,\dots,2d-2\}$, in the same way as we did before.
(This is the same notation as in Section~\ref{secCon}. In
Section~\ref{secGen}, this set was denoted by $\mathcal E_1$.) 
Figure~\ref{exemple} shows a possible $\mathcal R$-FPL
configuration inside $\mathcal R$ in the case $d=5$. 
There, the induced matching
is the matching $X_0$ from above, and the subset which encodes the 
vertical edges crossing
$\xi_1$ which are occupied by the $\mathcal R$-FPL configuration is the set
$\mathcal E_0$ from above. The right part of Figure~\ref{exemple}
shows the corresponding paths $\pi(X_0)$ and $\pi(\mathcal E_0)$. 
It should be noted that we do not suppose that $|\mathcal E|=d-1$; 
this will be a consequence of the subsequent considerations.

Given a path $\pi(.)$ (that is, either a path $\pi(X)$ or a path
$\pi(\mathcal E)$), we
let $\pi_i(.)$ be the $y$-coordinate of the point on the path which has
$x$-coordinate $i$. For example, we have $\pi_2(X_0)=1$
and $\pi_2(\mathcal E_0)=3$ (see Figure~\ref{exemple}).
Clearly, the two $y$-coordinates $\pi_i({\mathcal
E})$ and 
$\pi_i(X)$ differ by an even integer, so that we can define the integers
$$h_i=\frac{\pi_i({\mathcal E})-\pi_i(X)}{2},\quad i=1,2,\dots,2d-1.$$

\begin{figure}[htb]
\psfrag{O1}{$\mathcal{V}_1$}
\psfrag{O2}{$\mathcal{V}_2$}
\psfrag{O3}{$\mathcal{V}_3$}
\psfrag{O7}{$\mathcal{V}_7$}
\psfrag{O8}{$\mathcal{V}_8$}
\psfrag{O9}{$\mathcal{V}_9$}
\psfrag{O1}{$\mathcal{V}_1$}
\psfrag{C1}{$\varLambda_1$}
\psfrag{C2}{$\varLambda_2$}
\psfrag{C3}{$\varLambda_3$}         
\psfrag{C7}{$\varLambda_7$}
\psfrag{C8}{$\varLambda_8$}
\psfrag{C9}{$\varLambda_9$}
\psfrag{0}{$0$}
\psfrag{1}{$1$}
\centering \includegraphics[width=0.6\textwidth]{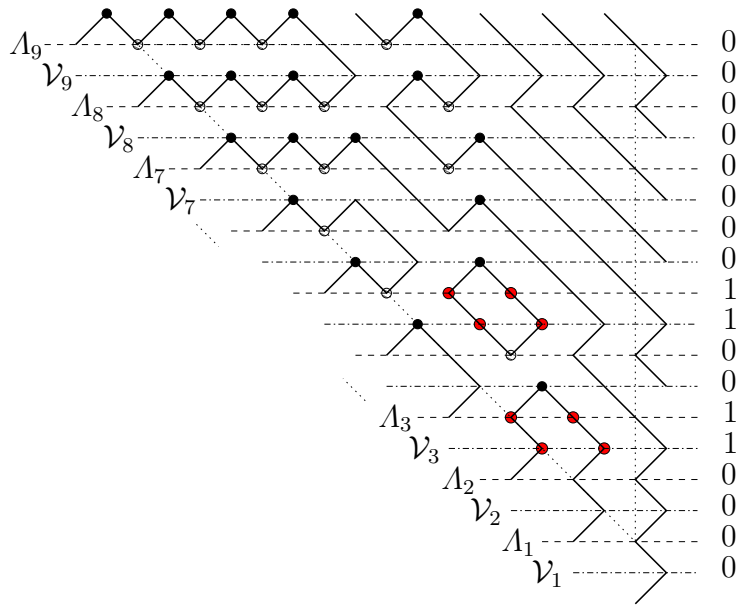}
\caption{The region $\mathcal{R}$ after rotation}
\label{Exemple_suite}
\end{figure}

In the sequel, for convenience,
we shall in fact consider a rotated picture, see
Figure~\ref{Exemple_suite} for the rotated version of the $\mathcal R$-FPL
configuration in $\mathcal R$ from Figure~\ref{exemple}. 
(The numbers on the right should be ignored at the moment.)
We always
also include the external edges occupied by the matching $X$.

\begin{figure}[htb]
\psfrag{c}{$(\land)$}
\psfrag{o}{$(\lor)$}         
\psfrag{s}{$(o)$}
\psfrag{a}{$(t)$}
\centering \includegraphics[width=0.8\textwidth]{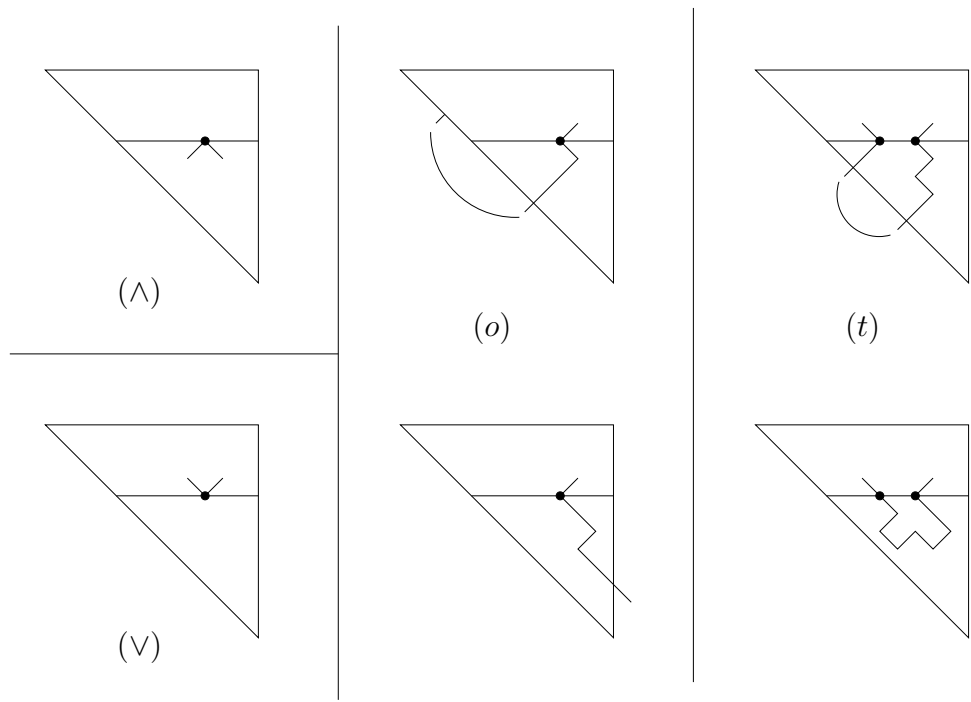}
\caption{The four types of vertices}
\label{types}
\end{figure}

We now cut the (rotated) region $\mathcal{R}$ in slices, 
as indicated in Figure~\ref{Exemple_suite}.
We call the cutting lines 
$\varLambda_i$ and $\mathcal{V}_i$, respectively, $i=1,2,\ldots, 2d-1$. 
Given an $\mathcal R$-FPL configuration,
we assign one of four {\it types} to every vertex on one of these lines,
as we explain below.
The reader is advised to look at Figure~\ref{types} while
reading the definition of these types.

Let $v$ be a vertex on line $l$, where $l$ is either $\varLambda_i$
or $\mathcal V_i$, for some $i$.
We say that $v$ has type $\lor $ (respectively $\land $) if
the two edges attached to $v$ are both above (respectively below) 
$l$.
If we are not in one of these cases, then there is exactly one edge
$e$ attached to $v$ which is below $l$. The edge $e$ is on a loop of the 
configuration. Let us start with $e$ and follow the part of the loop
below $l$. At some point, we shall either exit $\mathcal{R}$ or 
pass above $l$. Let us call the last edge before this happens $e'$.
The two
remaining types of vertices are defined depending on the location of
$e'$:
if $e'$ crosses $\xi_1$ (the right border of the
rotated version of $\mathcal R$), or is an external edge
linked, via $X$, to another external edge situated strictly above
$l$, then we say that $v$ has type $o$. (Here, ``$o$'' stands for
``out''.) If $e'$ is incident to $l$, or is an
external edge linked, via $X$, to another external edge situated
below $l$, then we say that $v$ has type $t$. 
See Figure~\ref{types}. An important feature of vertices
of type $t$ is that they always go in pairs. 
(This feature explains also the terminology
``$t$'', which stands for ``twin''.)
By convention, the right-most 
vertex of $\mathcal{V}_i$ is of type $o$ if it is not linked
to its neighbour on $\varLambda_{i-1}$.

Now we want to count the vertices on each line according to their
type: let $\Lambda_\lor (i)$, $\Lambda_\land (i)$, $\Lambda_t(i)$, and
$\Lambda_o(i)$ be the
respective numbers of vertices of $\varLambda_i$ of type $\lor $,
$\land $, $t$, and 
$o$. Likewise we introduce the quantities
$V_\lor (i)$, $V_\land (i)$, $V_t(i)$, and $V_o(i)$
for the line $\mathcal{V}_i$. Since vertices of type $t$ go in
pairs, $\Lambda_t(i)$ and $V_t(i)$ are even numbers. We are now ready to
state the result from which Lemmas~\ref{imposs}, \ref{lem:de} and
\ref{lem:3} will follow.

 \begin{thm} \label{th:main}
If $1\le i \le 2d-1$, we have
 $$\Lambda_\land (i)+\frac{1}{2}\Lambda_t(i)=V_\lor (i)+\frac{1}{2}V_t(i)=h_i.$$
\end{thm}

For an illustration of this result, the reader should consult
Figure~\ref{Exemple_suite}. The numbers along the right rim of 
the figure are the quantities $\Lambda_\land (i)+\frac{1}{2}\Lambda_t(i)$ and
$V_\lor (i)+\frac{1}{2}V_t(i)$. They are
indeed equal to the numbers
$h_i$, as can be seen by comparison with the right part 
of Figure~\ref{exemple}. 

Theorem~\ref{th:main} is an immediate corollary of the following two
lemmas, as can be easily seen by an induction with respect to $i$. In
the lemmas, we use the notations
$\DC(i)
 =(\Lambda_\land (i)+\frac{1}{2}\Lambda_t(i))- (V_\lor (i)+\frac{1}{2}V_t(i))$ and $\DO(i)
 =(V_\lor (i)+\frac{1}{2}V_t(i))- (\Lambda_\land
(i-1)+\frac{1}{2}\Lambda_t(i-1))$ for the differences of the
quantities which appear in the theorem.

\begin{lem} \label{lem1} For $i=1,2,\dots,2d-1$, we have $\DC(i)=0$.
\end{lem}

\begin{lem} \label{lem2} For $i=1,2,\dots,2d-1$, we have
$$\DO(i)= \begin{cases} \hphantom{-}1&\text{$s_i(\mathcal
E)$ is an up-step and $s_i(X)$ is a down-step,}\\
-1&\text{$s_i(\mathcal
E)$ is a down-step and $s_i(X)$ is an up-step,}\\
\hphantom{-}0&\text{otherwise,}\\
\end{cases}$$
where $s_i(\mathcal E)$ denotes the $i$-th step of the path
$\pi(\mathcal E)$, with the analogous meaning for $s_i(X)$.
\end{lem}


\begin{proof}[Proof of Lemmas \ref{lem1} and \ref{lem2}]
To begin with, we have to introduce some terminology and notation.
We consider the slices
between the lines $\mathcal{V}_i$ and $\varLambda_i$ on the one hand,
$i=1,2,\dots,2d-1$,
and between the lines $\varLambda_{i-1}$ and $\mathcal{V}_i$ on the
other hand, $i=2,3,\dots,2d-1$, separately. See Figure~\ref{straightened}.
One should note that the vertex on the very left of $\mathcal V_i$ is
always incident to an external edge, while the right-most vertices
of $\mathcal V_i$ and $\varLambda_i$ are connected by one of the
horizontal edges (before rotation) crossing $\xi_1$. These edges are
indicated as thick edges in Figure~\ref{straightened}.

\begin{figure}[htb]
\psfrag{Ci}{$\varLambda_i$} \psfrag{Oi}{$\mathcal{V}_i$}
\psfrag{Ci-1}{$\varLambda_{i-1}$} \centering
\includegraphics[width=0.8\textwidth]{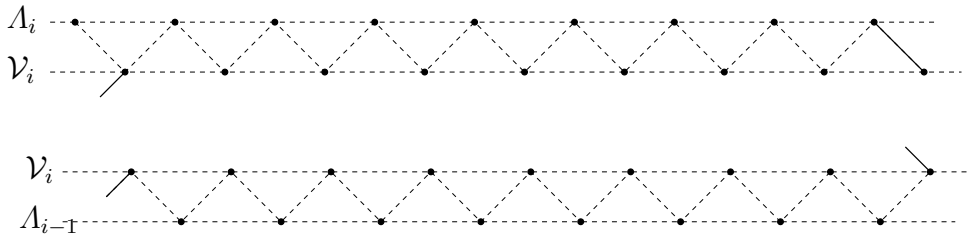}
\caption{The two kinds of slices}
\label{straightened}
\end{figure}

The restriction of an 
$\mathcal R$-FPL configuration to such a slice is a graph. Let
us call any connected component of such a graph which is not an
isolated vertex simply {\em
component}. Two components might be linked
in the sense that each of them has a vertex such that the two vertices
together form a pair of vertices of type $t$.
(The two vertices may be on the
same component, in which case this component is linked to
itself). When taking the reflexive and transitive closure of this
relation, we obtain an equivalence relation on the set of
components. We call the equivalence classes of this
relation {\em groups}. 
In other words, groups are minimal unions of components such that one
can pass from one component to another by a series of links.

For every subset of vertices $E$ of the graph, let us introduce the
notation $\DC_E$ (respectively $\DO_E$) 
to represent the quantity $\DC$ (respectively $\DO$) where
the count of vertices is restricted to those in $E$.  So, for $g$ a
group, $c$ a component, or for $Isol$ the set of isolated vertices, we
have the numbers $\DC_g$, $\DO_g$, $\DC_c$, $\DO_c$,
$\DC_{Isol}$, and $\DO_{Isol}$. As
all groups of a slice, together with the set of isolated vertices, 
form a partition of this slice, we have
$$
\DC=  \DC_{Isol} + \sum_{g \text{ a group}}\DC_g
\quad \text{and}\quad 
\DO=  \DO_{Isol} + \sum_{g \text{ a group}}\DO_g.
$$

In the sequel, when making pictures, we allow ourselves the following
simplifications. First, we write two linked components of a given
group next to each other even if they may be in a different order in
an actual configuration; some of the components may have to be
replaced by their left-right symmetric image in the process.  This is allowed
since this does not change the quantities $\Delta_g$. Second, we write
any component $c$ in the reduced way indicated in the column
``Simplification'' in Figure~\ref{comp}. Again, this is allowed since
this does not change the value of the corresponding $\Delta_c$.
An actual group and its simplified version are shown in
Figure~\ref{ex_simplif}.

\begin{figure}[htb]
\psfrag{I}{\small $I$}
\psfrag{II}{\small $II$}
\psfrag{III}{\small $III$}
\psfrag{comp}{\small Component}
\psfrag{simp}{\small Simplification}
\psfrag{coeff}{\small Coefficient}
\psfrag{0}{\small $\kern-5pt\hphantom{{}-}0$}
\psfrag{1}{\small $\kern-5pt\hphantom{{}-}1$}
\psfrag{-1}{\small $\kern-8pt{}-1$}                    
\centering
\includegraphics[width=0.6\textwidth]{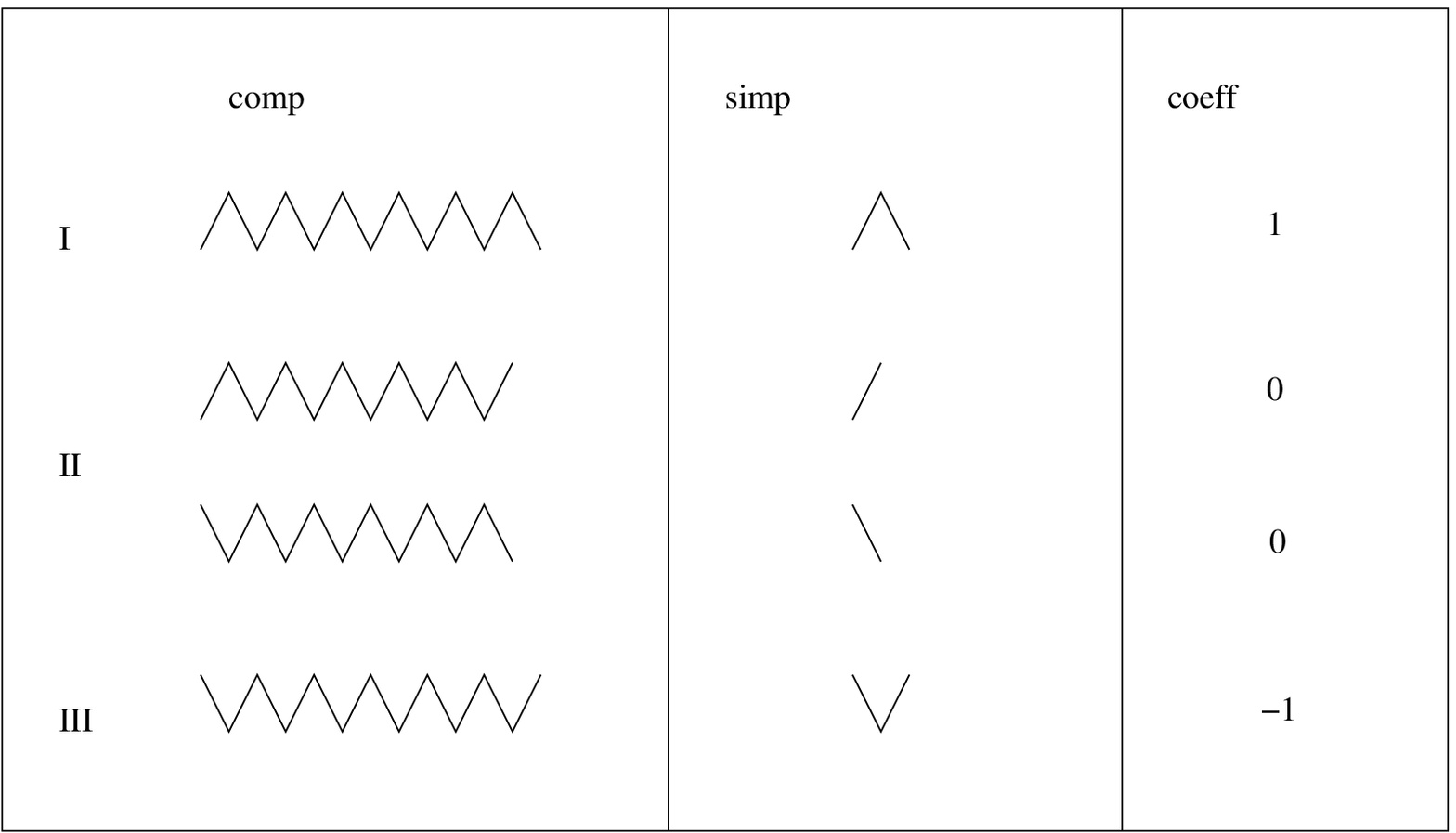}
\caption{The possible components, their simplified version, and the
coefficients used in the proof of Lemma~\ref{lem2}}
\label{comp}
\end{figure}

\begin{figure}[htb]
\psfrag{1}{$1$}
\psfrag{2}{$2$}
\psfrag{3}{$3$}
\psfrag{4}{$4$}
\centering
\includegraphics[width=0.6\textwidth]{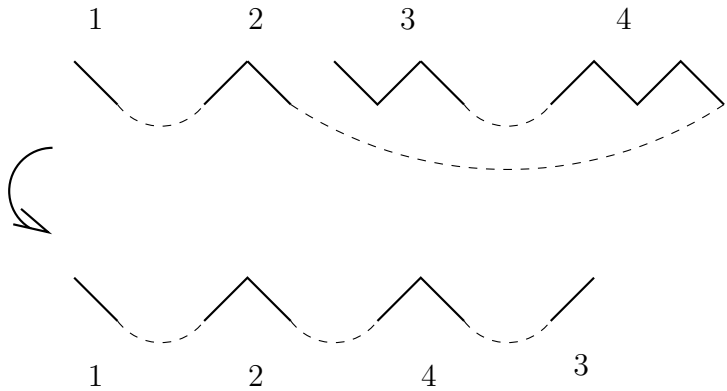}
\caption{A group (top) and its simplified version (bottom)}
\label{ex_simplif}
\end{figure}

\medskip
We are now ready for the proof of Lemma \ref{lem1}.
First, we claim that we have $\DC_{Isol}=0$, where $Isol$ 
denotes the set of isolated vertices in the slice
between $\mathcal{V}_i$ and $\varLambda_i$. 
Indeed, the isolated vertices in this slice are of type
$\land $ (respectively $\lor $) if they belong to $\mathcal{V}_i$
(respectively $\varLambda_i$). Thus,
they do not contribute to $\DC_{Isol}$, which is therefore $0$. 
Next we have to look
closer at groups in order to evaluate $\DC(i)$. The different
possibilities are illustrated in Figure~\ref{groups_1}. It is easy to
see that these are the only possible groups here, and that for each of
them we have $\DC_g=0$. Lemma~\ref{lem1} is thus proved.

\begin{figure}[htb]
\psfrag{k}{\small$(k)$}
\psfrag{k0}{\small$k\ge 0$}
\psfrag{k1}{\small$k\ge 1$}
\psfrag{x}{\small$x$}
\psfrag{s}{\small$o$}
\psfrag{y}{\small$y$}                     
\psfrag{wsf}{\small will stand for}
\psfrag{kcomp}{\small $k$ components}
\psfrag{a}{\small $t$}
\centering
\includegraphics[width=.5\textwidth]{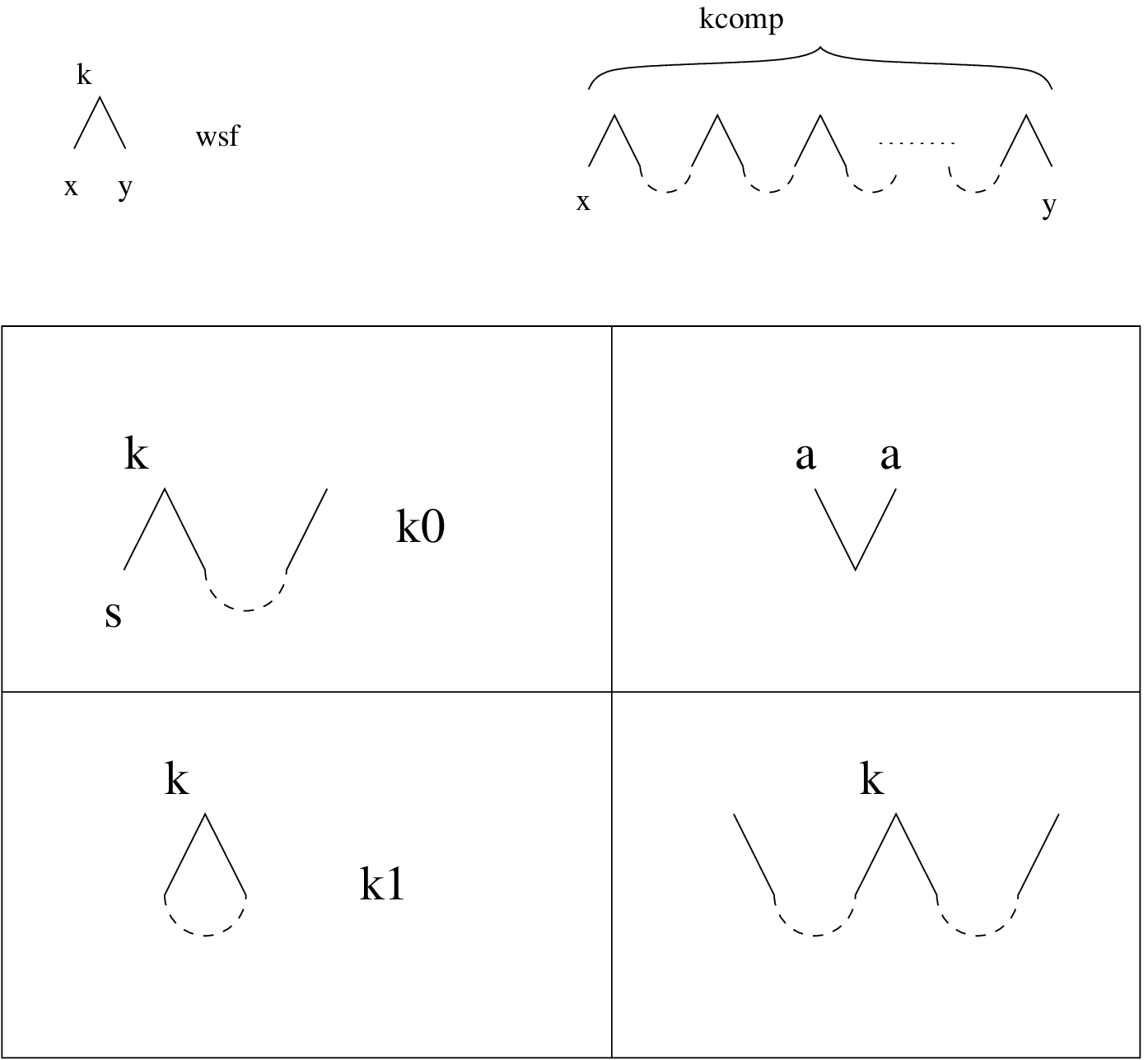}
\caption{All possible groups for Lemma \ref{lem1} } \label{groups_1}
\end{figure}

\medskip
Finally, we turn to the proof of Lemma \ref{lem2}, which is more
difficult, because for the isolated vertices $Isol$ in the slice
between  $\varLambda_{i-1}$ and $\mathcal{V}_i$ we do not have  
$\DO_{Isol}= 0$ in general, and
because there are more groups to consider.
First we use a trick that allows us to get rid of $\DO_{Isol}$ by
suitably modifying the  $\DO_g$'s.
Let us consider the slice between
$\varLambda_{i-1}$ and $\mathcal{V}_i$. It alternates between components and
blocks of isolated vertices. (When we say ``alternate,'' then we also
allow for empty blocks of isolated vertices.) A block $b$ of isolated vertices 
is generally surrounded
by two components. Notice that $\DO_b$ is equal to $1$
(respectively $-1$, respectively $0$) if the left-most vertex 
and the right-most vertex of $b$ are both
on $\mathcal{V}_i$ (respectively both on $\varLambda_{i-1}$,
respectively one on
each line). We may split these values as $1=\frac {1} {2}+\frac {1}
{2}$, $-1=(-\frac {1} {2})+(-\frac {1} {2})$, and $0=\frac {1} {2}+(-\frac
{1} {2})$, respectively $0=(-\frac {1} {2})+\frac {1} {2}$, and 
the idea is to ``move'' the left summand to the component to the left
and the right summand to the component to the right of $b$.
Thus, every component $c$ 
will accumulate one or two values out of $\{\frac {1} {2},-\frac
{1} {2}\}$. If the component $c$ is in between two blocks of isolated
vertices (recall: these blocks may also be empty blocks), then the 
sum of these two values is $1$, $0$, or $-1$. Column
``Coefficient'' in Figure~\ref{comp} shows this sum for the various
types of components. Thus, for such a component, we put
\begin{equation} \label{eq:Delta} 
\Delta'_c=\DO_c+\varepsilon,
\end{equation}
where $\varepsilon\in\{-1,0,1\}$ as indicated in the column ``Coefficient''
in Figure~\ref{comp}. If a component $c$ is to the very left of the
slice, then there is only a block of isolated vertices to the
right. In that case we also make the definition \eqref{eq:Delta}, but
with $\varepsilon=\frac {1} {2}$ if the right-most vertex of $c$ is on
$\varLambda_{i-1}$, and with 
$\varepsilon=-\frac {1} {2}$ if the right-most vertex of $c$ is on
$\mathcal{V}_i$. Similarly, if a component $c$ is to the very right of the
slice, then there is only a block of isolated vertices to the
left. In that case we also make the definition \eqref{eq:Delta}, but
with $\varepsilon=\frac {1} {2}$ if the left-most vertex of $c$ is on
$\varLambda_{i-1}$, and with 
$\varepsilon=-\frac {1} {2}$ if the left-most vertex of $c$ is on
$\mathcal{V}_i$. See Figure~\ref{groups_2} for these latter components.
(There, $s_i(\mathcal E)$ and $s_i(X)$ are the $i$-th steps of the
paths $\pi(\mathcal E)$ and $\pi(X)$, respectively, and ``up'' is short
for ``up-step,'', while ``down'' is short for ``down-step.'')
Finally, we put $\Delta'_g=\sum
\Delta_c'$. Then we have
\begin{equation} \label{eq:simpl}
\DO=\sum_{g\text{ a group}}\Delta'_g.
\end{equation}

\begin{figure}[htb]
\psfrag{delta}{\small$\Delta_g'$} \psfrag{grp}{\small Group $g$}
\psfrag{s}{\small $o$} \psfrag{0.5}{\small $\frac{1}{2}$}
\psfrag{-0.5}{\small $-\frac{1}{2}$} 
\psfrag{p1}{\small $\kern-5pts_i(X)=\text{up}$}
\psfrag{p0}{\small $\kern-7pt\hbox{$\begin{matrix} s_i(X)\\
      \quad =\text{down}\end{matrix}$}$} 
\psfrag{c1}{\small $s_i(\mathcal E)=\text{up}$} 
\psfrag{c0}{\small $\kern-6pt\hbox{$\begin{matrix} 
      s_i(\mathcal E)\\\quad =\text{down}\end{matrix}$}$} 
\psfrag{cp}{\small $\kern-25pts_i(X)=s_i(\mathcal E)=\text{down}$}
\psfrag{delt}{\small$\Delta_g'=0$} 
\centering \includegraphics[width=0.8\textwidth]{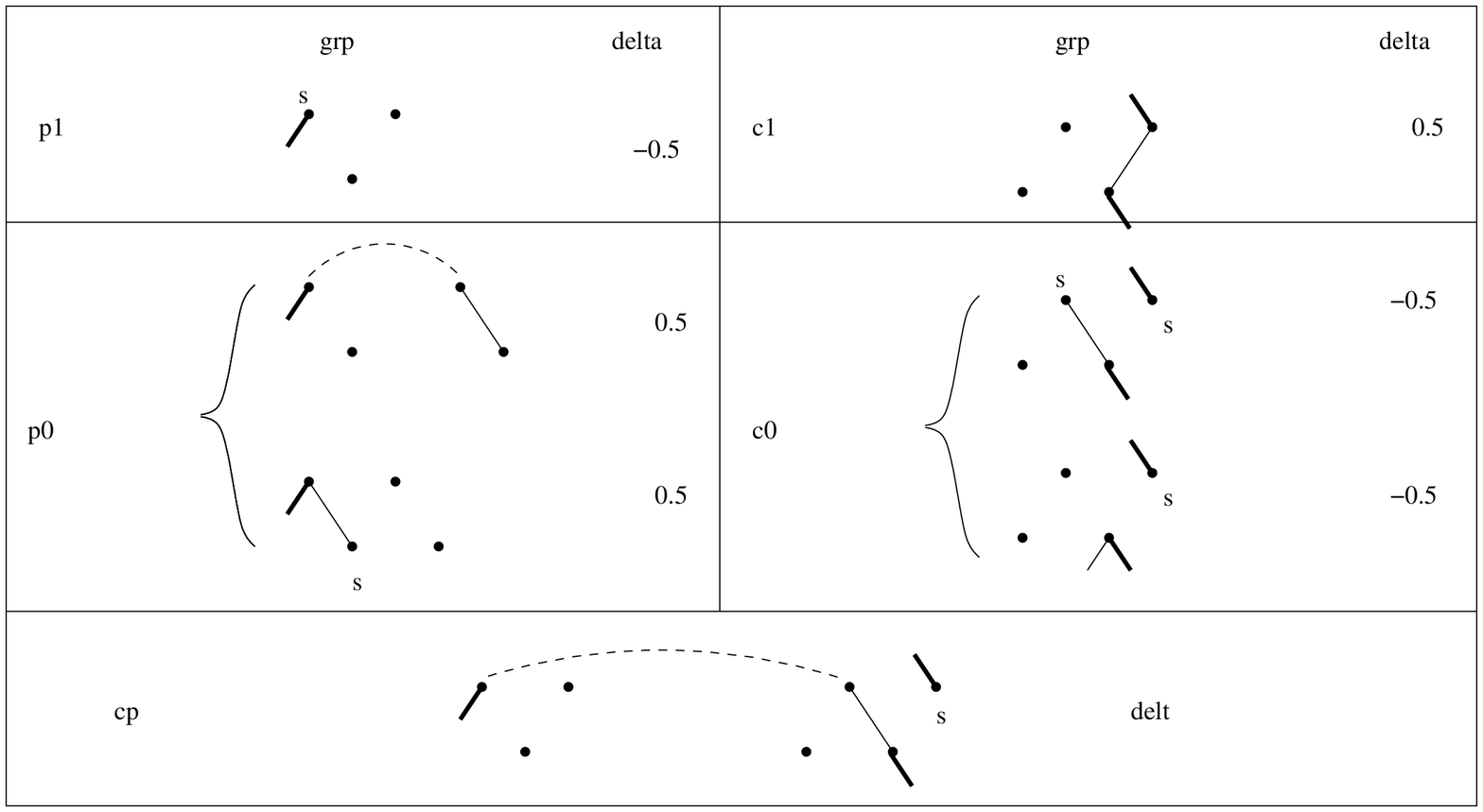}
\caption{``Side'' groups} \label{groups_2}
\end{figure}

Now for each group in Figure~\ref{groups_1} (that is, the groups
consisting exclusively of components 
between two blocks of isolated vertices), one verifies that
$\Delta'_g=0$. All possible cases for 
the remaining groups (that is, groups with components
to the very right or to the very left of the slice) are dealt with
in Figure~\ref{groups_2}. There, we use the
simplifications mentioned above plus an additional one: the pairs of
vertices of type $t$ 
on the line $\varLambda_{i-1}$ are not represented because they do
not modify the quantities $\Delta_g$.
It is not difficult to verify that, when the corresponding values of
$\Delta_g'$ given in 
Figure~\ref{groups_2} are substituted in (\ref{eq:simpl}), 
the assertion of Lemma~\ref{lem2} follows.
\end{proof}


\begin{cor} \label{cor}
\begin{enumerate} 
\item If there exists an $\mathcal R$-FPL configuration,
then the path $\pi(\mathcal E)$ is always weakly above the path
$\pi(X)$. Equivalently, $\la(\mathcal E)\subseteq \la(X)$.
\item If there exists an $\mathcal R$-FPL configuration,
then the path $\pi(\mathcal E)$ is a Dyck path. In particular,
$\mathcal E$ has exactly $d-1$ elements.
\item If $\pi(\mathcal E)=\pi(X)$ {\em(}or, equivalently, $\la(\mathcal
E)=\la(X)${\em)}, then there is exactly one 
$\mathcal R$-FPL configuration filling the
region $\mathcal R$. 
\end{enumerate}
\end{cor}

\begin{proof}
The assertion (1) is just a
reformulation of $h_i\ge 0$, which is a direct consequence of
Theorem~\ref{th:main}. 

For proving (2), we first observe that, because of (1), and because we
know that $\pi(X)$ is a Dyck path, $\pi(\mathcal
E)$ is certainly a path which never passes below the $x$-axis. Thus,
since the last steps of $\pi(X)$ and $\pi(\mathcal E)$ are down-steps
by definition, it suffices to prove that $h_{2d-1}=0$. To see the
latter, we look at the line $\varLambda_{2d-1}$. All its vertices are
incident to one of the horizontal edges along the right border of
$\mathcal R$ (the terminology refers to the
situation before rotation). Thus, none of these vertices can be of
type $\land$. Furthermore, since our $\mathcal R$-FPL configurations
contain either loops connecting external edges occupied by $X$ or
enter $\mathcal R$ from left and exit on the right, none of the
vertices on $\varLambda_{2d-1}$ can be of type $t$. Thus, by
Theorem~\ref{th:main}, we have $h_{2d-1}=\Lambda_\land
(i)+\frac{1}{2}\Lambda_t(i)=0$, as desired.

Finally, we prove assertion (3).
Let $\pi(\mathcal E)=\pi(X)$. In this case
we have $h_i=0$ for all $i$. Thus, according to Theorem~\ref{th:main},
the task is to construct $\mathcal R$-FPL configurations such that the lines
$\varLambda_i$ contain exclusively vertices of types $\lor $ and $o$,
while lines $\mathcal V_i$ contain exclusively vertices of types
$\land $ and $o$. We claim that,
if one proceeds slice by slice, beginning from the
bottom (that is, the slice between $\mathcal V_1$ and $\varLambda_1$),
and working upwards, then there is a unique such
configuration. 

To see this, the crucial observation is that, since
vertices of type $t$ are forbidden, there cannot be any closed loops
in $\mathcal R$ or loops which start and end on the top-most line
$\varLambda_{2d-1}$. Thus, inductively, if we have already worked our way up
to line $\mathcal V_i$, then, on $\mathcal V_i$, the
vertices of type $\land$ will be on the left and the vertices of type
$o$ will be on the right. (It is also possible that there are no
vertices of one of the two types.) 
See Figure~\ref{fig:V}. There, we have marked
the external edge incident to the left-most vertex on $\mathcal V_i$
and the horizontal edge (before rotation) crossing $\xi_1$ which
connects the right-most vertices of $\mathcal V_i$ and $\varLambda_i$
by thick lines. The latter edge forces us to connect the vertices of
type $o$ on $\mathcal V_i$ diagonally up-left to the neighbouring 
vertices of $\varLambda_i$. Since the latter cannot be of type
$\land$, they must be of type $o$. The remaining vertices on
$\varLambda_i$ (which are all on the left) must be of type $\lor$.

\begin{figure}[htb]
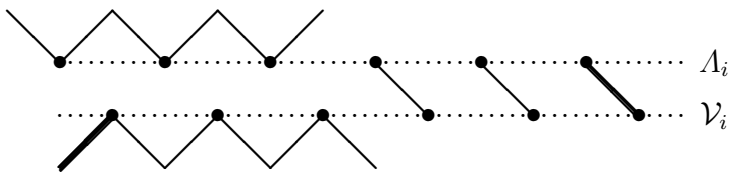

$$
\Einheit.7cm
\DuennPunkt(0,0)
\DuennPunkt(2,0)
\DuennPunkt(4,0)
\DuennPunkt(6,0)
\DuennPunkt(8,0)
\DuennPunkt(10,0)
\DuennPunkt(-1,1)
\DuennPunkt(1,1)
\DuennPunkt(3,1)
\DuennPunkt(5,1)
\DuennPunkt(7,1)
\DuennPunkt(9,1)
\thicklines
\Pfad(-1,-1),3\endPfad
\raise.8pt\hbox to0pt{\kern-.4pt\Pfad(-1,-1),3\endPfad\hss}%
\raise.4pt\hbox to0pt{\kern-.4pt\Pfad(-1,-1),3\endPfad\hss}%
\raise-.4pt\hbox to0pt{\kern0pt\Pfad(-1,-1),3\endPfad\hss}%
\raise-.8pt\hbox to0pt{\kern.4pt\Pfad(-1,-1),3\endPfad\hss}%
\Pfad(10,0),8\endPfad
\raise.4pt\hbox to0pt{\kern.4pt\Pfad(10,0),8\endPfad\hss}%
\raise-.8pt\hbox to0pt{\kern-.4pt\Pfad(10,0),8\endPfad\hss}%
\thinlines
\Pfad(0,0),43434\endPfad
\Pfad(-2,2),434343\endPfad
\Pfad(5,1),4\endPfad
\Pfad(7,1),4\endPfad
\SPfad(-1,0),111111111111\endSPfad
\SPfad(-1,1),111111111111\endSPfad
\Label\r{\mathcal V_i}(11,0)
\Label\r{\varLambda_i}(11,1)
\hskip7cm
$$
\caption{\label{fig:V}The unique filling of the slice between
$\mathcal V_i$ and $\varLambda_i$}
\end{figure}

On the other hand, again inductively,
if we have already worked our way up
to line $\varLambda_{i-1}$, then, on $\varLambda_{i-1}$, the
vertices of type $\lor$ will be on the left and the vertices of type
$o$ will be on the right. (It is also possible that there are no
vertices of one of the two types.) 
See Figures~\ref{fig:L1} and \ref{fig:L2}. There, we have marked
the external edge incident to the left-most vertex on $\mathcal V_{i}$
and the horizontal edges (before rotation) crossing $\xi_1$ which
are incident to 
the right-most vertices of $\varLambda_{i-1}$ and $\mathcal V_i$
by thick lines. Now there are two cases. Either the vertical edge
(before rotation) connecting the right-most vertices of
$\varLambda_{i-1}$ and $\mathcal V_i$ is among the vertical edges
chosen via $\mathcal E$ or not. In the first case
(cf.\ Figure~\ref{fig:L1}), this edge
forces us to connect the vertices of
type $o$ and as well the right-most vertex of type $\lor$ 
on $\varLambda_{i-1}$ diagonally up-right to the neighbouring 
vertices of $\mathcal V_i$. Since the latter vertices cannot be of type
$\lor$, they must be of type $o$. The remaining vertices on
$\mathcal V_i$ (which are all on the left) must be of type $\land$.
In the second case (cf.\ Figure~\ref{fig:L1}), the missing edge on the
very right 
forces us to connect the vertices of
type $o$ on $\varLambda_{i-1}$ diagonally up-left to the neighbouring 
vertices of $\mathcal V_i$. Since the latter vertices cannot be of type
$\lor$, they must be of type $o$, except for the left-most of those
which is of type $\land$. 
The remaining vertices on
$\mathcal V_i$ (which are all on the left) must also be of type $\land$.

\begin{figure}[htb]
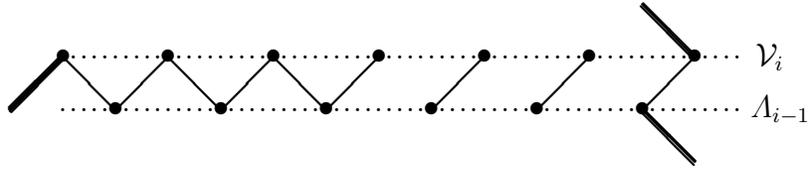

$$
\Einheit.7cm
\DuennPunkt(0,0)
\DuennPunkt(2,0)
\DuennPunkt(4,0)
\DuennPunkt(6,0)
\DuennPunkt(8,0)
\DuennPunkt(10,0)
\DuennPunkt(-1,1)
\DuennPunkt(1,1)
\DuennPunkt(3,1)
\DuennPunkt(5,1)
\DuennPunkt(7,1)
\DuennPunkt(9,1)
\DuennPunkt(11,1)
\thicklines
\Pfad(-2,0),3\endPfad
\raise.8pt\hbox to0pt{\kern-.4pt\Pfad(-2,0),3\endPfad\hss}%
\raise.4pt\hbox to0pt{\kern-.4pt\Pfad(-2,0),3\endPfad\hss}%
\raise-.4pt\hbox to0pt{\kern0pt\Pfad(-2,0),3\endPfad\hss}%
\raise-.8pt\hbox to0pt{\kern.4pt\Pfad(-2,0),3\endPfad\hss}%
\Pfad(11,1),8\endPfad
\raise.4pt\hbox to0pt{\kern.4pt\Pfad(11,1),8\endPfad\hss}%
\raise-.8pt\hbox to0pt{\kern-.4pt\Pfad(11,1),8\endPfad\hss}%
\Pfad(11,-1),8\endPfad
\raise.4pt\hbox to0pt{\kern.4pt\Pfad(11,-1),8\endPfad\hss}%
\raise-.8pt\hbox to0pt{\kern-.4pt\Pfad(11,-1),8\endPfad\hss}%
\thinlines
\Pfad(-1,1),434343\endPfad
\Pfad(6,0),3\endPfad
\Pfad(8,0),3\endPfad
\Pfad(10,0),3\endPfad
\SPfad(-1,0),1111111111111\endSPfad
\SPfad(-1,1),1111111111111\endSPfad
\Label\r{\mathcal V_i}(12,1)
\Label\r{\hbox{$\kern8pt\varLambda_{i-1}$}}(12,0)
\hskip7cm
$$
\caption{\label{fig:L1}The unique filling of the slice between
$\varLambda_{i-1}$ and $\mathcal V_i$ in Case 1}
\end{figure}

\begin{figure}[htb]


$$
\Einheit.7cm
\DuennPunkt(0,0)
\DuennPunkt(2,0)
\DuennPunkt(4,0)
\DuennPunkt(6,0)
\DuennPunkt(8,0)
\DuennPunkt(10,0)
\DuennPunkt(-1,1)
\DuennPunkt(1,1)
\DuennPunkt(3,1)
\DuennPunkt(5,1)
\DuennPunkt(7,1)
\DuennPunkt(9,1)
\DuennPunkt(11,1)
\thicklines
\Pfad(-2,0),3\endPfad
\raise.8pt\hbox to0pt{\kern-.4pt\Pfad(-2,0),3\endPfad\hss}%
\raise.4pt\hbox to0pt{\kern-.4pt\Pfad(-2,0),3\endPfad\hss}%
\raise-.4pt\hbox to0pt{\kern0pt\Pfad(-2,0),3\endPfad\hss}%
\raise-.8pt\hbox to0pt{\kern.4pt\Pfad(-2,0),3\endPfad\hss}%
\Pfad(11,1),8\endPfad
\raise.4pt\hbox to0pt{\kern.4pt\Pfad(11,1),8\endPfad\hss}%
\raise-.8pt\hbox to0pt{\kern-.4pt\Pfad(11,1),8\endPfad\hss}%
\Pfad(11,-1),8\endPfad
\raise.4pt\hbox to0pt{\kern.4pt\Pfad(11,-1),8\endPfad\hss}%
\raise-.8pt\hbox to0pt{\kern-.4pt\Pfad(11,-1),8\endPfad\hss}%
\thinlines
\Pfad(-1,1),43434\endPfad
\Pfad(6,0),8\endPfad
\Pfad(8,0),8\endPfad
\Pfad(10,0),8\endPfad
\SPfad(-1,0),1111111111111\endSPfad
\SPfad(-1,1),1111111111111\endSPfad
\Label\r{\mathcal V_i}(12,1)
\Label\r{\hbox{$\kern8pt\varLambda_{i-1}$}}(12,0)
\hskip7cm
$$
\caption{\label{fig:L2}The unique filling of the slice between
$\varLambda_{i-1}$ and $\mathcal V_i$ in Case 2}
\end{figure}

This completes the proof of the corollary.
\end{proof}

\end{document}